\documentclass[12pt]{article}
\usepackage{amsmath,amsthm,amssymb,amscd}
\usepackage{latexsym}
\newtheorem{thm}{Theorem}[section]
 \newtheorem{cor}[thm]{Corollary}
 \newtheorem{lem}[thm]{Lemma}
 \newtheorem{prop}[thm]{Proposition}
 \theoremstyle{definition}
 
 \theoremstyle{remark}

 \numberwithin{equation}{section}

\title
{Four-orbifolds with positive isotropic curvature}

\author{ Hong Huang}
\date{}
\begin{document}
\maketitle
\begin{abstract}
 We prove the following result: Let $(X,g_0)$ be a complete,
connected 4-manifold with uniformly positive isotropic curvature and
with bounded geometry. Then there is a finite collection
$\mathcal{F}$  of manifolds of the form  $\mathbb{S}^3 \times
\mathbb{R} /G$, where $G$ is a  discrete subgroup
of the isometry group of the round cylinder $\mathbb{S}^3\times
\mathbb{R}$ on which $G$ acts freely, such that $X$ is diffeomorphic to a possibly infinite
connected sum of  $\mathbb{S}^4,\mathbb{RP}^4$ and
members of $\mathcal{F}$. This extends  recent work of Chen-Tang-Zhu
and Huang. We also extend the above result to the case of orbifolds.
The proof uses Ricci flow with surgery on complete orbifolds.


\end{abstract}
\maketitle


\section {Introduction}

This is a continuation of our previous work [Hu1] on classifying open
4-manifolds with uniformly positive isotropic curvature.  Following Chen-Tang-Zhu's work [CTZ]  we'll  remove the condition of no essential
incompressible space form in [Hu1] and obtain the following

 \begin{thm} \label{thm 1.1}
  Let $(X,g_0)$ be a complete, connected 4-manifold with uniformly
positive isotropic curvature and with bounded geometry. Then there
is a finite collection $\mathcal{F}$ of manifolds of the form
$\mathbb{S}^3 \times \mathbb{R} /G$, where $G$ is a  discrete subgroup of the isometry group of the round cylinder $\mathbb{S}^3\times \mathbb{R}$ on which $G$ acts freely, such that $X$ is diffeomorphic
to a possibly infinite connected sum of
$\mathbb{S}^4,\mathbb{RP}^4$ and members of $\mathcal{F}$.
\end{thm}

 \noindent (By [MW] it is easy to see that the converse is also true:
 Any 4-manifold as in the conclusion of the theorem admits a complete metric
  with uniformly positive isotropic curvature and with bounded
geometry. The notion of  a possibly infinite connected sum will be
given later in this section.)

\noindent This also extends the Main Theorem in Chen-Tang-Zhu [CTZ] to
the noncompact case.

  As an  immediate consequence we have the following result which extends Corollary 2 in [CTZ].

\begin{cor} \label{cor 1.2}  \ \  A  4-manifold admits a complete, locally
conformally flat
 metric  with bounded geometry  and with
uniformly positive scalar curvature if and only if it admits a
complete metric with bounded geometry and with uniformly positive
isotropic curvature.
\end{cor}

Recall ([MM], [MW]) that a Riemannian manifold $M$ is said to have
uniformly positive isotropic curvature  if there exists a constant $c>0$ such that for all points $p \in M$ and all
orthonormal 4-frames $\{e_1,e_2,e_3,e_4\}\subset T_pM$, the curvature
tensor satisfies
\begin{equation*}R_{1313}+R_{1414}+R_{2323}+R_{2424}-2R_{1234}\geq c.
\end{equation*}
This notion can be easily adapted to the case of Riemannian
orbifolds.

Also recall that a complete Riemannian manifold  (or orbifold)  $M$ is
said to have bounded geometry if the sectional curvature is bounded
(in both sides) and the volume of any unit ball in $M$ is
uniformly bounded below away from zero.

Now we explain the notion of a possibly infinite connected sum which
slightly generalizes that in [BBM], [Hu1]. Let  $\mathcal{X}$ be a class
of smooth 4-manifolds. A smooth 4-manifold $X$ is said to be a possibly infinite
connected sum of members of $\mathcal{X}$ if there exists a
countable graph $G$ and a map $v\mapsto X_v$ which associates to
each vertex of $G$ a copy of some manifold in $\mathcal{X}$, such
that by removing from each $X_v$ as many open 4-balls as vertices
incident to $v$ and gluing the thus punctured $X_v$'s to each other (cf. pp. 102-106 of [BJ] or pp. 90-92 of [K])
according to the edges of $G$, one
obtains a smooth 4-manifold diffeomorphic to $X$. Note that we do
not assume that the elements in $\mathcal{X}$ are closed manifolds or the
graph is locally finite; compare [BBM], [Hu1]. By the way, note that by Cerf's theorem any diffeomorphism of the 3-sphere extends to the 4-ball.

 Inspired by Hamilton [H97], Perelman [P1], [P2], Chen-Zhu [CZ2], Chen-Tang-Zhu [CTZ], Bessi$\grave{e}$res et al [BBB$^+$], [BBM], and Kleiner-Lott [KL2], as in  [Hu1] we want to use a version of surgical Ricci flow  to prove Theorem 1.1. But as already noticed by
  [H97] (see also [CTZ]),  to do surgery on a  4-manifold with positive isotropic curvature
 (which contains essential incompressible space form) will lead in general to orbifolds (with isolated singularities).
 So we need to extend the construction in [Hu1] to the orbifold situation. Our treatment of surgery procedure is somewhat different from that in
 [CTZ] which considers the case of compact 4-orbifolds with at most isolated singularities. The main difference is that we do surgery slightly before (instead of exactly when) the curvature blows up. Of course, on the other hand, we also borrow many ideas and results from [CTZ].
 In fact we will prove the following  more general result than Theorem 1.1.

\begin{thm} \label{thm 1.3} \ \  Let $\mathcal{O}$ be a complete, connected
Riemannian 4-orbifold with uniformly positive isotropic curvature
and with bounded geometry. Then there is a finite collection
$\mathcal{F}$ of spherical 4-orbifolds  such that $\mathcal{O}$ is
diffeomorphic to a possibly infinite orbifold connected sum of
members of $\mathcal{F}$.
\end{thm}

When $\mathcal{O}$ is compact and with at most isolated
singularities the above theorem is due to [CTZ], see Theorem 2.1
there. The orbifold singularities in our Theorem 1.3 may not be isolated. We refer to Thurston [T] and Kleiner-Lott [KL2] for an
introduction to the topology and geometry of (effective) orbifolds. The notion
of orbifold connected sum will be explained in Section 2.  (For a 3-dimensional analogue of Theorem 1.3 see [Hu2], which generalizes some related results in [P2], [KL2] and [BBM].)

I expect that the construction of [MW] can be adapted to the case of orbifold connected sum to give  the converse of our Theorem 1.3:
 Any 4-orbifold as in the conclusion of  Theorem 1.3 should carry a complete metric
  with uniformly positive isotropic curvature and with bounded
geometry.   In fact one only needs to see whether any diffeomorphism of a spherical 3-orbifold $\mathbb{S}^3//\Gamma$ is isotopic to an isometry. In many cases this is indeed true (maybe this is already known in all cases, but I don't know):  By [Mc] in the spherical 3-manifold case this is true; by Theorem 2 of [CuZ], in the case  that $\Gamma <SO(4)$ and the exterior of the singular set admits a complete hyperbolic structure, this is also true.

   To consider Ricci flow with surgery on complete
4-orbifolds with not necessarily isolated singularities will encounter some difficulties which do not occur in the
 case of manifolds (see [H97], [P2], [CZ2], [BBM], and [Hu1]), or the case of compact orbifolds with at most isolated singularities (see [CTZ]), or the case of compact, orientable 3-orbifolds (see [KL2], where under normalized initial condition there is an a priori uniform upper bound on the orders of the isotropy groups for all time during the Ricci flow with surgery, see the discussion after Assumption 7.17 on p. 44 there, which is not the case here; moreover the relationship between the smooth category and the topological category is different in dimension 4, compared to that in  dimension 3). In [Hu1] we have established a crucial weak openness
(w.r.t. time) property of the canonical neighborhood condition for
the noncompact manifold case  (see Claim 1 in the proof of Proposition 3.6 there), which can be easily extended to the
noncompact orbifold case. We recall that the proof of Claim 1 in the proof of [Hu1, Proposition 3.6] uses the fact that since we do surgery when the scalar curvature reaches  some threshold (which is a function determined by the surgery parameters) before it blows up, the scalar curvature (hence the sectional curvature via Hamilton's pinching estimate [H97]) of the surgical solution is controlled (in terms of  the surgery parameters).

One of the main points in the present paper is to establish a canonical neighborhood structure for the restricted ancient $\kappa$--orbifold solutions in 4-dimension (see Proposition 3.6). For this purpose, we first
extend Gromoll-Meyer's theorem on complete,
noncompact manifold with positive
sectional curvature [GM] to the general orbifold case (see
Proposition 3.4); with the help of it, we use an argument involving  soul, distance function and Busemann function (which is somewhat different from that in [CTZ], the latter at some points only applies to the case of orbifolds with at most isolated singularities) to complete the proof of Proposition 3.6.

We will also generalize Theorem C2.4 of Hamilton
[H97] on gluing necks to the  orbifold situation (see
Proposition 2.3), which is very useful. (In the case of 4-orbifolds with at most isolated singularities, it suffices to use Hamilton's original version of gluing manifold necks.) For example, it is utilized to get Hamilton's canonical
uniformization for long tubes.
 With the aid of this uniformization we can pull back the orbifold Ricci flow solutions to manifolds, an idea already exploited in [CTZ] in the case of compact 4-orbifolds with isolated singularities. This will help us to overcome an additional difficulty
in the orbifold case, that is, the canonical neighborhoods in this case may a priori be very collapsed, cf. also [CTZ].  For another application of the gluing see also the proof of Proposition 4.4. We  also need to work slightly harder   to show that the noncollapsing property
survives the surgery  in the noncompact orbifold case than in the case of
compact 4-orbifolds with isolated singularities [CTZ] and in the case of noncompact 4-manifolds [Hu1]; see the proof of
Lemma 5.4. Meanwhile we will establish the property bounded curvature at bounded
distance (see Proposition 4.1) and persistence of almost standard
caps  (see Proposition 5.1) under the orbifold surgical Ricci flow, which are crucial
in the process of constructing $(r, \delta, \kappa)$-surgical
solutions (see Theorem 5.5).

In Section 2 we  introduce various notions on necks, in particular, Hamilton's canonical uniformization of necks, and prove a gluing result on necks.  In Section 3 we first prove the existence of Ricci flow on complete orbifolds with bounded curvature, then we describe the canonical neighborhood structure of
restricted ancient $\kappa$-solutions on 4-orbifolds. In Section 4 we choose the cutoff
parameters for surgical Ricci flow under the canonical neighborhood assumption, and describe
the metric surgery procedure. In Section 5, we construct  $(r, \delta,
\kappa)$-surgical solutions starting with a complete, connected
Riemannian 4-orbifold with uniformly positive isotropic curvature
and with bounded geometry. Finally, in Section 6, we prove
Theorem 1.3 and Theorem 1.1 using the construction in Section 5.  In most cases we will follow the notations and conventions in [BBB$^+$] and [Hu1].

\section{$\varepsilon$-necks and their gluing}

 To describe the structure of ancient $\kappa$-solution we need various notions on necks and caps. First we define topological necks and caps. Let $\Gamma$ be a finite
subgroup of $O(4)$. Following [CTZ], a (4-dimensional) topological neck is an orbifold which is
diffeomorphic to $\mathbb{S}^3//\Gamma \times \mathbb{R}$, where $\mathbb{S}^3//\Gamma $ denotes the quotient orbifold (as in [KL2]). Note that the notion topological neck has a slightly different meaning in Hamilton [H97] (compare Section C.2 in [H97], where the weaker condition local diffeomorphism is imposed). Let  $\sigma$ be an isometry of  $\mathbb{S}^3//\Gamma$ with $\sigma^2=1$, consider the quotient orbifold $(\mathbb{S}^3//\Gamma \times
\mathbb{R})//\{1,\hat{\sigma}\}$  (compare Remark 2.15 in [KL2]), where $\hat{\sigma}$ is the
involution on the orbifold $\mathbb{S}^3//\Gamma \times \mathbb{R}$
defined by $\hat{\sigma}(x,s)=(\sigma(x),-s)$ for $x\in
\mathbb{S}^3//\Gamma $ and $s\in \mathbb{R}$. Sometimes we also denote this
orbifold by $\mathbb{S}^3//\Gamma \times_{\mathbb{Z}_2}  \mathbb{R}$. (By the
way, note that we can consider $\Gamma$ and $\hat{\sigma}$ as isometries of $\mathbb{S}^4$ in a natural way, by lifting $\sigma$ to an isometry of $\mathbb{S}^3$ and viewing $\mathbb{S}^4$ as a suspension of
$\mathbb{S}^3$. We'll use the same notations for these isometries of $\mathbb{S}^4$.)  We define a topological cap to
be an orbifold diffeomorphic to either $\mathbb{S}^3//\Gamma
\times_{\mathbb{Z}_2}  \mathbb{R}$, or $\mathbb{R}^4//\Gamma$ with $\Gamma$  a finite subgroup of $O(4)$ (in this paper $\mathbb{R}^4$ is always endowed with the standard smooth structure unless explicitly stated otherwise; note that  it follows from [Hu1] that no exotic $\mathbb{R}^4$ can be endowed with a complete Riemannian metric with  uniformly
positive isotropic curvature and with bounded geometry).

We point out that our definitions of topological necks and caps are natural
extension of those in [CZ2] and [CTZ].
 For example, when $\Gamma$ is a finite subgroup of $O(4)$ which acts freely on $\mathbb{S}^3$  and such that $\mathbb{S}^3/\Gamma$ admits a fixed point free isometric involution
 $\sigma$, $(\mathbb{S}^3/\Gamma \times
\mathbb{R})//\{1,\hat{\sigma}\}$
defined as above
 is  the smooth cap
$C_{\Gamma}^\sigma$ defined in [CTZ]. For another example, let
$\sigma_i: \mathbb{S}^3 \rightarrow \mathbb{S}^3$ ($1\leq i \leq 4$)
be the four involutions defined by $(x_1,x_2,x_3,x_4) \mapsto
(x_1,-x_2,-x_3,-x_4)$, $(x_1,x_2,x_3,x_4) \mapsto
(x_1,x_2,-x_3,-x_4)$,  $(x_1,x_2,x_3,x_4) \mapsto
(x_1,x_2,x_3,-x_4)$ and  $(x_1,x_2,x_3,x_4) \mapsto
(x_1,x_2,x_3,x_4)$ respectively. Then the quotients $(\mathbb{S}^3
\times \mathbb{R})//\{1,\hat{\sigma_i}\}$  defined as above are
orbifolds with singular set of dimension 0, 1, 2 and 3 respectively.
Note that $(\mathbb{S}^3 \times \mathbb{R})//\{1,\hat{\sigma_1}\}$
is the orbifold cap of type II defined in [CTZ], denoted by $\mathbb{S}^4/(x,\pm
  x')\setminus \bar{\mathbb{B}}^4$ there. Also note that when
$\mathbb{R}^4//\Gamma$ has only  an isolated singularity, it is the same as the orbifold cap of type I
defined in [CTZ], denoted by $C_\Gamma$ there.

It turns  out that any (topological) neck or cap can be written as  a so called
infinite orbifold connected sum of spherical 4-orbifolds. Our
notion of possibly infinite orbifold connected sum
 extends both the manifold case (as
defined in Section 1), the finite orbifold connected sum defined in
[CTZ] (where the orbifolds considered have at most isolated
singularities), and the operation of performing 0-surgery defined in
[KL2]. Let $\mathcal{O}_i$ ($i=1,2$) be two $n$-orbifolds, and let
$D_i \subset \mathcal{O}_i$ be two embedded suborbifolds-with boundary, both
diffeomorphic to some quotient orbifold $D^n//\Gamma$, where $D^n$
is the closed unit $n$-ball, and $\Gamma$ is a finite subgroup of
$O(n)$. Choose a diffeomorphism $f:
\partial D_1 \rightarrow \partial D_2$, and use it to glue together
$\mathcal{O}_1\setminus int (D_1 )$ and  $\mathcal{O}_2\setminus int
(D_2 )$. The result is called the orbifold connected sum of
$\mathcal{O}_1$ and $\mathcal{O}_2$ via the gluing map $f$, and is
denoted by $\mathcal{O}_1 \sharp_f  \mathcal{O}_2 $. If $ D_i$
($i=1,2$) are disjoint embedded suborbifolds-with boundary (both
diffeomorphic to some quotient orbifold $D^n//\Gamma$) in the same connected
$n$-orbifold $\mathcal{O}$, the result of similar process as above
is called the orbifold connected sum on (the single orbifold)
$\mathcal{O}$, and is denoted by $\mathcal{O}\sharp_f$.

 Given a collection $\mathcal{F}$ of $n$-orbifolds, we say an $n$-orbifold $\mathcal{O}$ is a possibly infinite
 orbifold connected sum of  members of $\mathcal{F}$  if there exist a  countable  graph $G$ (in which we allow an edge to connect some vertex
to itself), a map $v\mapsto F_v$ which associates to each vertex of
$G$ a copy of some orbifold in $\mathcal{F}$, and a map $e \mapsto
f_e$ which associates to each edge of $G$ a self-diffeomorphism of
some $(n-1)$-dimensional spherical orbifold, such that if we do an
orbifold connected sum of the corresponding $F_v$('s)  via the
gluing map $f_e$ for each edge $e$, we obtain an $n$-orbifold diffeomorphic to
$\mathcal{O}$.

By the way, note that  the graph $G$ which describes  the possibly
infinite orbifold connected sum appearing in Theorem 1.3 is locally
finite, since the members of collection $\mathcal{F}$ in this
theorem are compact.

 Using a somewhat ambiguous notation,  we can write the necks and caps  as
infinite orbifold connected sums of spherical 4-orbifolds as mentioned above:

 $\mathbb{S}^3//\Gamma \times \mathbb{R}^1 \hspace*{1mm} \approx \hspace*{1mm}
\cdot\cdot\cdot \hspace*{1mm} \sharp \hspace*{1mm} \mathbb{S}^4//\Gamma  \hspace*{1mm}  \sharp \hspace*{1mm}
\mathbb{S}^4//\Gamma \hspace*{1mm} \sharp \hspace*{1mm} \cdot\cdot\cdot$,

$\mathbb{R}^4//\Gamma \hspace*{1mm} \approx \hspace*{1mm}
\mathbb{S}^4
//\Gamma \hspace*{1mm} \sharp \hspace*{1mm} \mathbb{S}^4 //\Gamma \hspace*{1mm} \sharp \hspace*{1mm}  \cdot\cdot\cdot$, \hspace*{4mm} and

$(\mathbb{S}^3//\Gamma \times
\mathbb{R})//\{1,\hat{\sigma}\} \hspace*{1mm} \approx \hspace*{1mm}
\mathbb{S}^4//\langle \Gamma,\hat{\sigma}\rangle \hspace*{1mm}  \sharp \hspace*{1mm} \mathbb{S}^4//\Gamma
\hspace*{1mm} \sharp \hspace*{1mm} \mathbb{S}^4//\Gamma \hspace*{1mm} \sharp \hspace*{1mm} \cdot\cdot\cdot$.

Note that the
orbifold connected sums appearing in these  three examples are
actually the operation of performing 0-surgery as defined in [KL2]
(which can be extended to the non-oriented case), and we have
omitted the $f$'s in the notation.  Also note that for a
diffeomorphism $f: \mathbb{S}^3//\Gamma \rightarrow
\mathbb{S}^3//\Gamma$  the mapping torus $\mathbb{S}^3//\Gamma
\times_f \mathbb{S}^1 \approx \mathbb{S}^4//\Gamma \sharp_f$. By
work of [CuZ] and [Mc], in the case $\Gamma < SO(4)$, the mapping
class group of $\mathbb{S}^3//\Gamma$ is finite. So given a finite
subgroup $\Gamma < SO(4)$, there are only a finite number of
orbifolds of the form $\mathbb{S}^3//\Gamma \times_f \mathbb{S}^1$
up to diffeomorphism.  Note that
every orientation-reversing
homeomorphism of $\mathbb{S}^3$ has a fixed point,  so if $\Gamma < O(4)$ is a finite subgroup such that  $\mathbb{S}^3//\Gamma$ is a manifold, then
$\Gamma < SO(4)$.
\hspace *{0.4cm}

Now following Perelman, we define $\varepsilon$-neck, $\varepsilon$-cap, and strong
$\varepsilon$-neck. As in Definition 2.20 in [KL2], we do not
require the map in the definition of $\varepsilon$-closeness of two
pointed orbifolds to be precisely basepoint-preserving.  Given a
Riemannian 4-orbifold $(\mathcal{O},g)$, an open subset $U$, and a
point $x_0\in U$. $U$ is an $\varepsilon$-neck centered at $x_0$ if
there is a diffeomorphism $\psi: (\mathbb{S}^3//\Gamma)
  \times \mathbb{I} \rightarrow U$  such that the pulled back metric
  $\psi^*g$, scaling with some factor $Q$, is $\varepsilon$-close (in
  $C^{[\varepsilon^{-1}]}$ topology) to the standard metric on $(\mathbb{S}^3//\Gamma)
  \times \mathbb{I}$ with scalar curvature 1 and
  $\mathbb{I}=(-\varepsilon^{-1},\varepsilon^{-1})$, and the distance $d(x_0,  |\psi|(\mathbb{S}^3//\Gamma \times \{0\}))< \varepsilon/ \sqrt{Q}$.
   (Here $\Gamma$ is a finite subgroup of isometries of $\mathbb{S}^3$.)  By the way, note that the notion $\varepsilon$-neck here is somewhat stricter than the notion geometrically $(\varepsilon, k)$ cylindrical neck in Hamilton [H97].
    An open subset
  $U$ is an $\varepsilon$-cap centered at $x_0$ if  $U$ is diffeomorphic to $\mathbb{R}^4//\Gamma$ or  $\mathbb{S}^3//\Gamma
\times_{\mathbb{Z}_2} \mathbb{R}$, and there is an open set $V$ with
compact closure such that $x_0 \in V \subset \overline{V}\subset U$,
and $U\setminus \overline{V}$ is an $\varepsilon$-neck.
              Given a 4-dimensional orbifold Ricci flow  $(\mathcal{O}, g(t))$, an open subset
$U$, and a point $x_0\in U$.
   $U$ is a strong
  $\varepsilon$-neck  centered at $(x_0,t_0)$ if    there
  is a diffeomorphim $\psi: (\mathbb{S}^3//\Gamma)
  \times \mathbb{I} \rightarrow U$ such that, the pulled back
  solution $\psi^*g(\cdot,\cdot)$ on the parabolic region $\{(x,t)| x \in U, t\in [t_0-Q^{-1},t_0]\}$ (for some $Q>0$),
  parabolically rescaled with factor $Q$, is $\varepsilon$-close (in $C^{[\varepsilon^{-1}]}$
  topology) to the subset $(\mathbb{S}^3//\Gamma
  \times \mathbb{I})\times [-1,0]$ of the evolving round cylinder $\mathbb{S}^3//\Gamma
  \times \mathbb{R}$, with scalar curvature one and length 2$\varepsilon^{-1}$ to $\mathbb{I}$ at time zero, and
  the distance at time $t_0$, $d_{t_0}(x_0,  |\psi|(\mathbb{S}^3//\Gamma \times \{0\}))<\varepsilon/ \sqrt{Q}$.

\hspace *{0.4cm}

Following Hamilton [H97], we introduce another refined notion on necks which is very useful for us. Let $(\mathcal{O}, g)$ be a Riemannian 4-orbifold and $\Phi: \mathbb{S}^{3}\times
(a,b)\rightarrow \mathcal{O}$  be a smooth map whose image  is a
suborbifold diffeomorphic to $\mathbb{S}^{3}//\Gamma
\times (a,b)$ for some finite subgroup $\Gamma <O(4)$. Suppose $\Phi=\bar{\Phi}\circ \pi_1$, where $\pi_1: \mathbb{S}^{3}\times
(a,b)\rightarrow  \mathbb{S}^{3}//\Gamma \times (a,b)$ is the natural projection,   and $\bar{\Phi}:\mathbb{S}^{3}//\Gamma
\times (a,b) \rightarrow \mathcal{O} $ is a  topological neck. Following [H97], we define the mean radius $r(z)$ of the horizontal sphere $\mathbb{S}^3 \times \{z\}$ in the pull-back metric $\Phi^*g$ (restricted to the sphere) so that its area (w.r.t. $\Phi^*g$) is $A(\mathbb{S}^3   \times \{z\}, \Phi^*g)=\sigma_3r(z)^3$, where $\sigma_3$ is the area of the unit round $\mathbb{S}^3$. We'll call $\Phi$  Hamilton's canonical uniformization  if in addition it satisfies the
conditions (a), (b), (c) and (d) listed in Section C.2  of [H97]. For the convenience of the readers, below we will quote these conditions from [H97].

\hspace *{0.1cm}

(a) every horizontal sphere $\mathbb{S}^3 \times \{z\}$ for $z\in (a,b)$ has constant mean curvature in the pull-back metric $\Phi^*g$;

(b) the identity map from every horizontal sphere in the standard  metric $\bar{g}$ on the cylinder  (restricted to the sphere) to the same sphere in the   pull-back metric $\Phi^*g$ (also restricted to the sphere) is harmonic;

(c) the volume of any subcylinder in the pull-back metric $\Phi^*g$ is given by
\begin{displaymath}
vol(\mathbb{S}^3   \times [w,w'], \Phi^*g)=\sigma_3\int_w^{w'}r(z)^4dz
\end{displaymath}
(compare Remark 3.10 (iii) in [HS]);  and

 (d) if $\bar{V}$ is a Killing vector field on $\mathbb{S}^3 \times \{z\}$ in the metric $\bar{g}$ restricted to the sphere, then
\begin{displaymath}
\int_{\mathbb{S}^3 \times \{z\}}\bar{g}(\bar{V},W)d\bar{a}=0
\end{displaymath}
for any unit vector field $W$ which is $\Phi^*g$-orthonormal to the sphere, where $d\bar{a}$ is the surface measure on the sphere induced by $\bar{g}$.

\hspace *{0.1cm}

Let $\Phi$ be a Hamilton's canonical uniformization with $\Phi=\bar{\Phi}\circ \pi_1$ as above, following [CTZ], the map $\bar{\Phi}$ is called  Hamilton's canonical parametrization.   Since as noted above our notion of topological neck is slightly different from that in Hamilton [H97], the notion Hamilton's canonical uniformization is also slightly different from the corresponding notion normal neck in [H97]. In particular, Hamilton allows (the image of the part near) the two ends of his normal neck to  overlap each other, while we do not. Anyway, the distinction between the two notions is very small, except that here we are dealing with the larger orbifold category. Virtually all theorems (with their proofs) in Section C.2 of
[H97] extend to our situation with minor modifications.

The following uniqueness result shows the rigidity of Hamilton's canonical uniformization.

\begin{lem} \label{lem 2.1}  There exists $\hat{\varepsilon}$ with the following property. If $\Phi_i$ ($i=1,2$) are two Hamilton's canonical  uniformizations in $(\mathcal{O},g)$ with $\Phi_i^*g$ locally $\hat{\varepsilon}$-close to the standard metrics on the cylinders after suitably rescaled, and with $\Phi_2=\Phi_1\circ F$, where $F$ is a diffeomorphism of the cylinders, then $F$ is an isometry in the standard metrics on the cylinders.
\end{lem}

\noindent {\bf  Proof} \ \  It's trivial to adapt the proof of Lemma C2.1 in [H97] to our situation.
\hfill{$\Box$}

\hspace *{0.4cm}

We also have a result on the existence of Hamilton's canonical uniformization.

\begin{lem} \label{lem 2.2}  Given $\hat{\varepsilon} >0$, there exists $\varepsilon > 0$ such that if $U$ is an $\varepsilon$-neck centered at $x_0$ in a Riemannian 4-orbifold $(\mathcal{O},g)$ with a diffeomorphism   $\psi: (\mathbb{S}^3//\Gamma)
  \times (-\varepsilon^{-1}, \varepsilon^{-1}) \rightarrow U$, then   there exists Hamilton's canonical uniformization  $\Phi: \mathbb{S}^3 \times (-l,l)\rightarrow
U$, whose image contains
  the portion $\psi(\mathbb{S}^3 //\Gamma
\times (-0.98\varepsilon^{-1}, 0.98\varepsilon^{-1}))$ in $U$;  moreover, $\Phi^*g$ is $\hat{\varepsilon}$-close to the standard metric on the cylinder after suitably rescaled.
\end{lem}

\noindent {\bf  Proof} \ \ One can easily adapt the arguments in the proof of Theorem C2.2 and Corollary C2.3 in [H97] to our situation (compare  Lemma A.1 in [Hu1]). \hfill{$\Box$}

\hspace *{0.4cm}

Then we have the following proposition on gluing necks which extends Theorem
C2.4 in [H97], and which will be used several times later.

\begin{prop} \label{prop 2.3}
  There exist  $\varepsilon >0$  with the
  following property. Let $U_i$ ($i=1,2$) be two $\varepsilon$-necks in a 4-orbifold $(\mathcal{O}, g)$ with diffeomorphisms   $\psi_i: (\mathbb{S}^3//\Gamma_i)
  \times (-\varepsilon^{-1}, \varepsilon^{-1}) \rightarrow U_i$. Suppose one end of the neck $U_1$ does not intersect  $U_2$. Let $\Phi_i: \mathbb{S}^3 \times (-l_i,l_i)\rightarrow
U_i$ be Hamilton's canonical uniformization as above, whose image contains
  the portion $\psi_i(\mathbb{S}^3 //\Gamma_i
\times (-0.98\varepsilon^{-1}, 0.98\varepsilon^{-1}))$ in $U_i$ respectively. Suppose there is some point $P_1$ in the domain
  cylinder of $\Phi_1$ at standard distance at least $0.01\varepsilon^{-1}$ from the
  ends such that the point $\Phi_1(P_1)$ is also in the image of $\Phi_2$.
  Then  there exists Hamilton's canonical uniformization  $\Phi$ of the union $
\mbox{image} (\Phi_1)\cup \mbox{image} (\Phi_2 )$, and  diffeomorphisms $F_1$ and $F_2$
  of the cylinders such that $\Phi_1=\Phi\circ F_1$ and $\Phi_2=\Phi\circ F_2$.

\end{prop}

\noindent {\bf  Proof} \ \  We adapt the argument of Theorem C2.4 in
[H97] to our orbifold case, with the aid of the orbifold covering theory in
Section 2.2 of [BMP]. By Lemma 2.2  we may assume that $\Phi_i^*g$ is $\hat{\varepsilon}$-close to the standard metrics on the cylinders after suitably rescaled, where $\hat{\varepsilon}>0$ is small depending on $\varepsilon$. Let $P_2$ be a point in the domain cylinder of
$\Phi_2$ such that $\Phi_2(P_2)=\Phi_1(P_1)$. We may assume that the point $P:=\Phi_1(P_1)$ is not a
orbifold singular point, since otherwise we can choose an ordinary point near
$P$ which is also in the intersection of the image of $\Phi_1$ with the image of
$\Phi_2$. Assume that $P_i$ lies on the sphere $\mathbb{S}^3 \times
\{z_i\}$, $i=1,2$. We claim that we can find a map
$$G: \mathbb{S}^3 \times \{z_2\} \rightarrow \mathbb{S}^3 \times
\{z_1\}$$ \noindent  so that $\Phi_1\circ G=\Phi_2$ and $G(P_2)=P_1$. The reason
is as follows. Given any point
 $Q_2 \in \mathbb{S}^3 \times \{z_2\}$, choose any smooth path $\gamma_2$ in $ \mathbb{S}^3 \times \{z_2\}$
 from $P_2$ to $Q_2$ such that the image of $\gamma:=\Phi_2\circ \gamma_2$  contains at most finite number of orbifold singular points.
  Then by Section 2.2 of [BMP] we can lift $\gamma$ to a unique
 smooth path $\gamma_1$ (starting from $P_1$) in the first cylinder with $\Phi_1  \circ \gamma_1=\gamma$, and  reach a point $Q_1$ with $\Phi_1(Q_1)=\Phi_2(Q_2)$, since $\gamma_1$ is almost horizontal (as is easily seen), and if the path $\gamma_2$ is not too long and $\varepsilon$ (hence $\hat{\varepsilon}$) is sufficiently small, $\gamma_1$ cannot run out of the first cylinder (since its length is nearly the same as that of $\gamma_2$).
 (Note that we require the  lift $\gamma_1$ to be smooth, otherwise it may not be unique; also note that by our assumption  $P_1$ is (a little)  far away from the ends of the first cylinder.) By Section 2.2 of [BMP] (see the first paragraph on p. 36 in [BMP], which can be  proved by adapting the standard covering theory in, for example, Lemma 3.3 in Chapter 5 of [M]), the point $Q_1$ is independent of the choice of path $\gamma_2$ (satisfying the same
 condition), and the map $G: \mathbb{S}^3 \times \{z_2\} \rightarrow \mathbb{S}^3 \times
      (-l_1, l_1)$ taking  $Q_2$ to $Q_1$ is well-defined. By construction  we have $\Phi_1\circ G=\Phi_2$, and locally $G$ extends to an isometry from $\Phi_2^*g$ to $\Phi_1^*g$, so $G(\mathbb{S}^3 \times \{z_2\})$ is a constant mean curvature sphere in the first cylinder with the metric $\Phi_1^*g$. Furthermore the constant mean curvature sphere $G(\mathbb{S}^3 \times \{z_2\})$ is nearly horizontal and  passes through $P_1$, and must be the same as $\mathbb{S}^3 \times \{z_1\}$.

 Now  we can flip one of the cylinders $\mathbb{S}^3 \times (-l_i,l_i)$ (that is, reverse the $z$-direction) if necessary, and  extend $G$ in the same way as above  mapping $\mathbb{S}^3 \times [z_2,z_2+\mu]$ to $\mathbb{S}^3 \times [z_1,z_1+\mu]$ (or mapping $\mathbb{S}^3 \times [z_2-\mu,z_2]$ to $\mathbb{S}^3 \times [z_1-\mu,z_1]$) for some $\mu>0$. By Lemma 2.1, if $\varepsilon$ (hence $\hat{\varepsilon}$) is sufficiently small, $G$ is an isometry in the standard metrics of the cylinders. Then we can use $G$ to glue the cylinders  $\mathbb{S}^3 \times (-l_i,l_i)$ together, and get the desired $\Phi$ and $F_i$ from $\Phi_i$. \hfill{$\Box$}

\hspace *{0.4cm}

\noindent {\bf Remark} It is interesting to compare Proposition 2.3 and its proof to the statements and arguments in the Appendix of [MT] and in Section 3.2 of [BB$B^+$] for gluing smooth 3-dimensional necks. After suitably modifying the definition of Hamilton's canonical uniformization (in particular, adding the condition $(b')$ in Section C.2 of [H97]) one can obtain an analogue of Proposition 2.3 in the case of 3-orbifolds, which was used in [Hu2].

\hspace *{0.4cm}

\begin{prop} \label{prop 2.4}\ \ Let $\varepsilon$ be sufficiently small.
Let $(\mathcal{O},g)$ be a complete, connected Riemannian 4-orbifold. If each
point of $\mathcal{O}$ is the center of an $\varepsilon$-neck or an
$\varepsilon$-cap, then $\mathcal{O}$ is diffeomorphic to
a neck, a cap, or an orbifold connected
sum of at most two spherical orbifolds.
\end{prop}

\noindent {\bf  Proof}. Using Proposition
2.3, arguing as in the proof of Proposition 2.6 in [Hu1] we see that if each
point of $\mathcal{O}$ is the center of an $\varepsilon$-neck or an
$\varepsilon$-cap, then $\mathcal{O}$ is diffeomorphic to a neck, a cap, a mapping torus $\mathbb{S}^3//\Gamma
\times_f \mathbb{S}^1$,  or a union of two caps along their ends. As observed before, a mapping torus $\mathbb{S}^3//\Gamma \times_f \mathbb{S}^1$ is diffeomorphic to an orbifold connected sum on a single spherical 4-orbifold. Using the description of caps as orbifold connected sums of spherical 4-orbifolds given before, we see that a union of two caps along their ends is diffeomorphic to an orbifold connected sum of two spherical 4-orbifolds:

$\mathbb{R}^4//\Gamma \hspace*{1mm} {\cup}_f \hspace*{1mm}  \mathbb{R}^4 //\Gamma \hspace*{1mm} \approx \hspace*{1mm} \mathbb{S}^4 //\Gamma \hspace*{1mm} \sharp \hspace*{1mm} \mathbb{S}^4 //\Gamma$,

$(\mathbb{S}^3//\Gamma \times \mathbb{R})//\{1,\hat{\sigma}\} \hspace*{1mm} {\cup}_f \hspace*{1mm} \mathbb{R}^4//{\Gamma} \hspace*{1mm} \approx \hspace*{1mm}
\mathbb{S}^4//\langle \Gamma,\hat{\sigma}\rangle \hspace*{1mm} \sharp \hspace*{1mm} \mathbb{S}^4 //\Gamma$, and

$(\mathbb{S}^3//\Gamma \times \mathbb{R})//\{1,\hat{\sigma}\} \hspace*{1mm} {\cup}_f \hspace*{1mm} (\mathbb{S}^3//\Gamma \times \mathbb{R})//\{1,\widehat{\sigma'}\} \hspace*{1mm} \approx \hspace*{1mm}
\mathbb{S}^4// \langle \Gamma,\hat{\sigma}\rangle \hspace*{1mm} \sharp \hspace*{1mm} \mathbb{S}^4// \langle\Gamma,\widehat{\sigma'}\rangle$.

 \hfill{$\Box$}

\hspace *{0.4cm}

 It is clear
that there is a function $\varepsilon \mapsto f_3(\varepsilon)$ with
$f_3(\varepsilon) \rightarrow 0$ as $\varepsilon \rightarrow 0$, such
that if $N$ is an $\varepsilon$-neck with scaling factor $\lambda$
and $x, y \in N$, then one has

\begin{equation*}
|\lambda^{-1}R(x)-1|\leq f_3(\varepsilon), \hspace*{8mm}
|\frac{R(x)}{R(y)}-1| \leq f_3(\varepsilon).
\end{equation*}
(For example, see [BBB$^+$] Section 3.2, our notation  $f_3(\varepsilon)$ is
borrowed from there.)

  Let $K_{st}$ be the superemum of the sectional curvatures of the
(4-dimensional ) smooth standard solution on $[0,4/3]$. The
following lemma on strengthening necks extends [Hu1, Lemma A.2]; compare [BBB$^+$, Lemma
4.3.5] and [BBM, Lemma 4.11]; we emphasize that the condition (i) here is slightly weaker than the corresponding ones in these cited references, and is more flexible.
It will be needed later to fix our constants.

\begin{lem} \label{lem 2.5}
  For any $\varepsilon \in (0, 10^{-4})$
 there exists $\beta=\beta (\varepsilon)\in (0,1)$ with the
 following property.

  Let $a, b$ be real numbers satisfying $a< b <0$ and $|b| \leq
  \frac{3}{4}$, let $(\mathcal{O}(\cdot), g(\cdot))$ be a surgical solution to the Ricci flow (this notion will be recalled in Section 4)
  defined on $(a,0]$, and $x$ be a point such that:

  (i) $|R(x,b)-1|\leq f_3(\beta \varepsilon)$;

  (ii) $(x,b)$ is the center of a strong $\beta \varepsilon$-neck;

  (iii) $P(x, b, (\beta \varepsilon)^{-1}, |b|)$ is unscathed and
  satisfies $|Rm|\leq 2K_{st}$.

 \noindent  Then $(x,0)$ is the center of a strong $\varepsilon$-neck.
\end{lem}

\noindent {\bf Proof}  We  argue by contradiction. Otherwise there exist
$\varepsilon \in (0, 10^{-4})$, a sequence $\beta_k \rightarrow 0$,
sequences $a_k < b_k$, $b_k \in [-3/4,0]$, and a sequence of
surgical solution $(\mathcal{O}_k(t), g_k(t))$ ($t \in (a_k,0]$) with a point
$x_k \in \mathcal{O}_k$ such that

(i) $|R(x_k,b_k)-1|\leq f_3(\beta_k \varepsilon)$;

  (ii) $(x_k,b_k)$ is the center of a strong $\beta_k \varepsilon$-neck $N_k$ with a diffeomorphism $\psi_k: \mathbb{S}^3//\Gamma_k \times
(-(\beta_k \varepsilon)^{-1},(\beta_k \varepsilon)^{-1})\rightarrow  N_k$;

  (iii) $P(x_k, b_k, (\beta_k \varepsilon)^{-1}, |b_k|)$ is unscathed and
  satisfies $|Rm|\leq 2K_{st}$, but

  (iv) $(x_k,0)$ is not the center of any strong $\varepsilon$-neck.

\noindent For each $k$, let $\tilde{\psi_k}: \mathbb{S}^3 \times (-(\beta_k \varepsilon)^{-1},(\beta_k \varepsilon)^{-1})
\rightarrow  N_k$ be the composition of  the natural projection from $\mathbb{S}^{3}\times
(-(\beta_k \varepsilon)^{-1},(\beta_k \varepsilon)^{-1})$ to $ \mathbb{S}^{3}//\Gamma_k \times (-(\beta_k \varepsilon)^{-1},(\beta_k \varepsilon)^{-1})$  with $\psi_k$. Then we pull back
$(\mathcal{O}_k(t), g_k(t))$ ($t \in (a_k,0]$) to $\mathbb{S}^3 \times
(-(\beta_k \varepsilon)^{-1},(\beta_k \varepsilon)^{-1})$  via $\tilde{\psi_k}$. Now we can proceed as in the proof of
 [BBB$^+$, Lemma 4.3.5].   \hfill{$\Box$}

\section{ Restricted ancient $\kappa$-solutions on 4-orbifolds}

 Let $(\mathcal{O},g_0)$ be a complete orbifold
 with $|Rm|\leq K$.
Consider the Ricci flow ([H82])

\begin{displaymath}\frac{\partial g}{\partial t}=-2Ric, \ \  g|_{t=0}=g_0. \ \ \ \
 \end{displaymath}
Shi's short time existence for Ricci flow with initial data
a complete (noncompact) manifold with bounded sectional curvature ([S]) extends
to the orbifold case. This should be well-known to the experts, but since we cannot find a proof in the literature, we will indicate how to adapt Shi's original proof to the orbifold case for convenience.

\begin{thm} \label{thm 3.1}  Let $(\mathcal{O},g_0)$ be a complete $n$-orbifold
 with $|Rm|\leq K$.  Then the Ricci flow with initial data $(\mathcal{O},g_0)$
 has a short time solution with bounded curvature.
\end{thm}

\noindent {\bf Proof}  The compact case is well-known, see for example [H03]. We only need to consider the noncompact case.  Following [S], the idea is to convert the original problem to  solving a sequence of Dirichlet boundary value problems for Ricci-DeTurck flow ([D]) on exhausting domains. (For a nice exposition of Ricci-DeTurck flow see Chapter 3 of [CK], which we would like to follow when we adapt Sections 2 and 6 of [S] to the orbifold case.) We will only
indicate the necessary modifications to Shi [S]. First one can write the complete
$n$-orbifold $\mathcal{O}$ as a union of an increasing sequence of compact
$n$-suborbifolds with boundary embedded $(n-1)$-suborbifolds. (This is
possible, as is seen for example by the following argument. First we use a standard technique to construct a smooth proper function $f:=\sum_{j=1}^{\infty}j\eta_j$ on $\mathcal{O}$, where $\eta_j$ is a sequence of cut-off functions  adapted to some suitably chosen open cover $\{U_j\}$ of $\mathcal{O}$, then we apply Sard theorem and preimage theorem (see [BB]) to $f$ to get the desired result.) Note that the Hessian comparison theorem used in Section 4 of [S] holds true for orbifolds. (Compare Borzellino and Zhu [BZ].) Also note that Stokes theorem holds (so that we can integrate by parts) in a
bounded domain (with boundary an embedded $(n-1)$-suborbifold) in a Riemannian
orbifold. (In fact it holds for a slightly more general domain, see
[C]. ) Note that to carry out the kind of integration by parts as in Section 6 of [S] one only needs a corollary of Stokes theorem which also applies to the nonorientable case; compare for example, Theorem 14.34 (or Theorem 16.48 in the second edition) of [L]. In Lemma 3.1 (on p. 244) of [S], the dependence on the
injectivity radius of $\overline{D}$ can be replaced by that on the Sobolev
constant in some Sobolev inequality which holds in a bounded domain
(with boundary an embedded $(n-1)$-suborbifold) in a Riemannian orbifold; compare [LSU]. (For Sobolev
inequalities on compact manifolds with boundary, one can see [A] and
[He]; for Sobolev inequalities  on closed orbifold, see [Ch], [Na] and [F].
The extension to the case of compact orbifolds with boundary is
routine.) On pages 260 and 286  of [S], one can pull back the solution to
(a suitable ball in) $T_{\tilde{x}_0}\widetilde{U}$  via
$exp_{x_0}\circ \pi_*$, where $(\widetilde{U}, G, \pi)$ is a
 uniformizing chart for some $U \ni x_0$.
 \hfill{$\Box$

 \vspace *{0.4cm}

By extending the proof in [Ko] of Chen-Zhu's uniqueness theorem for Ricci flow on complete manifolds [CZ1] to the orbifold case,
the solution is unique in the category of bounded curvature
solutions  (even in a slightly larger category).

Now we  restrict to the 4-dimensional case. Let $(\mathcal{O},g_0)$ be a Riemannian 4-orbifold $(\mathcal{O},g_0)$ with
uniformly positive isotropic curvature.
  If we decompose the bundle $\Lambda^2T\mathcal{O}$
into the direct sum of its self-dual and anti-self-dual parts
\begin{equation*}
\Lambda^2T\mathcal{O}=\Lambda_+^2T\mathcal{O} \oplus
\Lambda_-^2T\mathcal{O},
\end{equation*}
 then the curvature operator can be decomposed as
\begin{equation*}
\mathcal{R}=\left(
  \begin{array}{cc}
    A & B \\
    B^T & C \\
  \end{array}
\right).
\end{equation*}
  Denote the eigenvalues of the matrices
$A, C$ and $\sqrt{BB^T}$ by $a_1\leq a_2 \leq a_3$, $c_1 \leq c_2
\leq c_3$ and $b_1 \leq b_2 \leq b_3$ respectively. It is easy to
see (cf. Hamilton [H97]) that  for a Riemannian 4-manifold/ orbifold
the condition of uniformly positive isotropic curvature is
equivalent to that
 there is a positive constant $c$ such that $a_1+a_2\geq c$,
$c_1+c_2\geq c$ everywhere.

We can easily generalize Hamilton's pinching result in [H97] to our
situation.

  \begin{lem} \label{lem 3.2}  \ \  (cf. Hamilton [H97])\ \ Let  $(\mathcal{O},g_0)$ be a complete 4-orbifold  with
  uniformly positive isotropic curvature  ($(a_1+a_2)\geq c$, $(c_1+c_2)\geq c$) and with bounded curvature ($|Rm|\leq K$).  Then there exist positive
  constants $\varrho, \Psi, L, P, S<+\infty$ depending only on the initial
  metric  (through $c, K$), such that any complete solution to the Ricci flow with initial data $(\mathcal{O},g_0)$ and with bounded
  curvature satisfies
 \begin{equation}
 \begin{aligned}
& a_1+\varrho>0, \hspace{8mm} c_1+\varrho>0, \\
 & \max\{a_3,b_3,c_3\}\leq \Psi(a_1+\varrho), \hspace{8mm}  \max\{a_3,b_3,c_3\}\leq
  \Psi(c_1+\varrho), \\
& \frac{b_3}{\sqrt{(a_1+\varrho)(c_1+\varrho)}}\leq 1+\frac{L
e^{Pt}}{\max \{\ln \sqrt{(a_1+\varrho)(c_1+\varrho)},S\}}
\end{aligned}
\end{equation}
\noindent at all points and times.
\end{lem}

Since the 4-orbifolds we consider have uniformly positive isotropic
curvature, and in particular, have uniformly positive scalar
curvature, any  Ricci flow staring with them  will blow up in finite time. It
follows from Lemma 3.2 that any blow-up limit (if it exists)
satisfies the following restricted isotropic curvature pinching
condition
$$a_3\leq \Lambda a_1, \hspace{4mm} c_3\leq \Lambda c_1,\hspace{4mm}  b_3^2\leq a_1c_1.$$

This motivates the following definition.

\hspace *{0.4cm}

{\bf Definition} (compare [CZ2], [CTZ]) Let $(\mathcal{O},g(t))$ be a
smooth, complete, and nonflat solution to the Ricci flow on a
4-orbifold.  It is said to be a (restricted) ancient
$\kappa$-orbifold solution if the following holds:

(i) it exists on the time interval $(-\infty, 0]$,

(ii) it has positive isotropic curvature and bounded sectional
curvature, and satisfies

$$a_3\leq \Lambda a_1, \hspace{4mm} c_3\leq \Lambda c_1, \hspace{4mm} b_3^2\leq a_1c_1$$

\noindent for some constant $\Lambda>0$, and

(iii) it is $\kappa$-noncollapsed on all scales for some constant
$\kappa>0$.

\vspace *{0.4cm} Now we will investigate  the structure of
(restricted) ancient $\kappa$-orbifold solutions.

\begin{prop} \label{prop 3.3}  \ \   Let $(\mathcal{O},g(t))$ be a
(restricted) ancient $\kappa$-orbifold solution such that the
curvature operator has null eigenvector somewhere. Then
$\mathcal{O}$ is isometric to a shrinking Ricci soliton
$\mathbb{S}^3//\Gamma \times \mathbb{R}$ or $\mathbb{S}^3//\Gamma
\times_{\mathbb{Z}_2} \mathbb{R}$ for some finite subgroup $\Gamma$  of $O(4)$.
\end{prop}

 \noindent {\bf Proof}\ \   Compare the proof of Theorem 3.4 in [CTZ]. We pull back our solution  to the
universal cover and get $(\widetilde{\mathcal{O}},\tilde{g}(t))$. Using
Hamilton's strong maximum principle ([H86], which can be adapted to
the orbifold case) and an orbifold de Rham decomposition theorem (see [KL2, Lemma 2.19]) we see that at any time $\widetilde{\mathcal{O}}$ is
isometric to a product $\mathcal{O}' \times \mathbb{R}$. Using the
definition of ancient $\kappa$-orbifold solution we see that at any time
$\mathcal{O}'$ is a round 3-sphere (compare  [CZ2, Lemma 3.2]). Now we see that $(\mathcal{O},g(t))$ is
a metric quotient of the evolving round cylinder $\mathbb{S}^3 \times \mathbb{R}$; compare for example, [R, Theorem 13.3.10].
 Then it follows from the $\kappa$-noncollapsing assumption that
$\mathcal{O}$ is noncompact (cf. [CZ2, p.212] and [CTZ, p. 52]). So $\mathcal{O}$ has one or
two ends. If $\mathcal{O}$ has two ends, it must be isometric to
$\mathbb{S}^3//\Gamma \times \mathbb{R}$ for some finite subgroup
$\Gamma < O(4)$. If $\mathcal{O}$ has one (and only one) end, it must be isometric
to $\mathbb{S}^3//\Gamma \times_{\mathbb{Z}_2} \mathbb{R}$ for some finite
subgroup $\Gamma < O(4)$. The reason (for the latter case) is as
follows.  We can write
$\mathcal{O}=\mathbb{S}^3 \times \mathbb{R}//\widetilde{\Gamma}$ for a subgroup
$\widetilde{\Gamma}$ of isometries of the round cylinder $\mathbb{S}^3 \times
\mathbb{R}$. Since $\mathcal{O}$ has one end, we can write $\widetilde{\Gamma}=\Gamma \cup \Gamma^1$, where the second
components of $\Gamma$ and $\Gamma^1$ act on $\mathbb{R}$ as the
identity or a reflection respectively. Since $\mathcal{O}$ has  only one end, $\Gamma^1 \neq \emptyset$. Pick $\sigma \in \Gamma^1$.
Then $\sigma^2 \in \Gamma$, and $\sigma \Gamma=\Gamma^1$. It follows that $\sigma$ induces an involution, denoted by $\bar{\sigma}$, acting isometrically on $\mathbb{S}^3//\Gamma \times \mathbb{R}$.   Now we see that
$\mathcal{O}=(\mathbb{S}^3//\Gamma \times \mathbb{R})//\langle
\bar{\sigma} \rangle$, which is clearly of the form $\mathbb{S}^3//\Gamma\times_{\mathbb{Z}_2} \mathbb{R}$.

 \hfill{$\Box$

\vspace*{0.4cm}

The following proposition extends Gromoll-Meyer's theorem [GM] and
Corollary 3.8 in [CTZ].

\begin{prop} \label{prop 3.4}  \ \  Let $(\mathcal{O},g)$ be a complete,
noncompact and connected Riemannian $n$-orbifold with positive
sectional curvature, then a soul of $\mathcal{O}$ is a point, and
$\mathcal{O}$ is diffeomorphic to $\mathbb{R}^n//\Gamma$ for some
finite subgroup $\Gamma$ of $O(n)$.
\end{prop}

\noindent {\bf Proof} \ \  Recall that a compact, totally convex suborbifold
$\mathcal{S}$ of $\mathcal{O}$ without boundary is called a soul. We
also recall the construction of soul in [CE], [CG], [GW] and [KL2].  Fix a base point $\star
\in |\mathcal{O}|$, and associate the
function $b_\eta(p):=\lim_{t\rightarrow +\infty}(d(p,\eta(t))-t)$ to any unit-speed ray $\eta$ in
$|\mathcal{O}|$ starting from $\star$.
Let the Busemann
function $F=\inf_{\eta}b_\eta$, where $\eta$ runs over unit speed rays
starting from $\star$. Then $F$ is a proper concave function on
$|\mathcal{O}|$ which is bounded above. The  superlevel sets of $F$
are compact totally convex subsets of $\mathcal{O}$. Now let $C_0$
be a compact totally convex subset  of $\mathcal{O}$ with $\partial
C_0 \neq \emptyset$. For example, one can choose $C_0$ as some
superlevel set of $F$. Inductively, for $i\geq 0$, let $C_{i+1}$ be
the minimal nonempty superlevel set of $d_{\partial C_i}:=d(\cdot,
\partial C_i)$ on $C_i$ if $\partial C_i \neq
\emptyset$. (For definition of $\partial C_i$ see Section 3 of
[KL2].) Let $S$ be the nonempty $C_k$ so that $C_{k+1}$ does not
exist (i.e. $\partial C_k =
\emptyset$), and let  $\mathcal{S}=\mathcal{O}|_S$. (Such $k$ exists and of course $k \ge 1$.) Then $\mathcal{S}$ is
a soul of $\mathcal{O}$, and $\mathcal{O}$ is diffeomorphic to the
normal bundle $\mathcal{NS}$ of $\mathcal{S}$ by Proposition 3.13 in
[KL2].

We need only to show that a soul  $\mathcal{S}$ as constructed above
is a point. We argue by contradiction. Otherwise, we can choose a
minimizing geodesic $\gamma:[a,b]\rightarrow S$. Since $S=C_k$ is
the minimal nonempty superlevel set of  $d_{\partial C_{k-1}}$ on
$C_{k-1}$, $S=d_{\partial C_{k-1}}^{-1}[c,+\infty)=d_{\partial
C_{k-1}}^{-1}\{c\}$ for some nonnegative constant $c$.  It follows
that $d_{\partial C_{k-1}}\circ \gamma=c$. Note that $c\neq 0$ (that
is $\gamma ([a,b])  \nsubseteq
\partial C_{k-1}$),
since $C_{k-1}\setminus \partial C_{k-1} \neq \emptyset$.

Let $t\mapsto exp_{\gamma(a)}tX(a)$ be a minimizing unit-speed
geodesic from $\gamma(a)$ to $\partial C_{k-1}$, defined for $t\in [0
,c]$. Let $\{X(s)\}_{s\in
[a,b]}$ be the parallel transport of $X(a)$ along $\gamma$. Then the
rectangle $V: [a,b]\times [0,c]\rightarrow C_{k-1}$ given by
$V(s,t)=exp_{\gamma(s)}tX(s)$ is flat and totally geodesic by Lemma
3.9 in [KL2]. This  contradicts  the positive sectional curvature
assumption.  \hfill{$\Box$

 \vspace*{0.4cm}

The following proposition extends Proposition 3.7 in [CTZ].

\begin{prop} \label{prop 3.5} \ \ There exist a universal positive constant $\eta$ and  a universal positive function $\omega: [0,\infty)\rightarrow (0,\infty)$
such that for any restricted ancient $\kappa$-orbifold solution
$(\mathcal{O},g(t))$, we have

\begin{displaymath}
(i) \hspace{4mm}  R(x,t)\leq R(y,t)\omega(R(y,t)d_t(x,y)^2)
\end{displaymath}
\noindent for any $x,y\in \mathcal{O}$, $t\in (-\infty,0]$, and

\begin{displaymath}(ii) \hspace{4mm}  |\nabla R|(x,t)\leq \eta R^{\frac{3}{2}}(x,t), \hspace{4mm}  |\frac{\partial R}{\partial
t}|(x,t)< \eta R^2(x,t)
\end{displaymath}
\noindent for any $x\in \mathcal{O}$, $t\in (-\infty,0]$.

 \end{prop}

\noindent {\bf Proof} \ \ (i) We follow the proof of Proposition 3.7 in
[CTZ].

Case 1: The curvature operator has null eigenvector somewhere. Then
the result follows immediately from Proposition 3.3.

Case 2: $\mathcal{O}$ is compact and has positive curvature
operator. In this case we know that $\mathcal{O}$ is diffeomorphic
to a spherical 4-orbifold by Hamilton [H86]  (see also B$\ddot{o}$hm-Wilking [BW] for a more general result under a weaker condition, where they also considered the orbifold case). (One can
also give a proof of this fact by combining  Perelman's
noncollapsing theorem [P1] and the Ricci flow compactness theorem
for orbifolds as established in [KL2] with Hamilton's estimates in
[H86].) Then we can pull back the solution to $\mathbb{S}^4$ and use
the corresponding result in the manifold case as established in
[CZ2], and argue as in [CTZ] to get the desired result.

Case 3: $\mathcal{O}$ is noncompact and has positive curvature
operator. By  Proposition 3.4 we know that $\mathcal{O}$ is
diffeomorphic to $\mathbb{R}^4//\Gamma$ for some finite subgroup
$\Gamma$ of $O(4)$. Then we argue as in Case 2.

(ii)  We argue by pulling back
the solutions to their universal covers (which are manifolds as
already shown) and using the corresponding result in the manifold case
in [CZ2].  \hfill{$\Box$

\hspace *{0.4cm}

In the proof of the following Proposition 3.6, we'll use a simple fact about the soul of a complete (noncompact) Riemannian manifold of positive sectional curvature with symmetry. The fact is implicit in the proof of Theorem 3 in [GM]. For the convenience of the readers, we'll provide a proof essentially following [GM].

 \hspace *{0.4cm}

 \noindent {\bf Fact} Let $M$ be a complete noncompact Riemannian manifold of positive sectional curvature, and $G$ be a group of isometries of $M$. Then there exists a point in $M$ which is a soul of $M$ and which is fixed by the group $G$.

 \hspace *{0.4cm}

 \noindent {\bf Proof of Fact} In the proof of  Theorem 3 in [GM] it is shown that there exists a compact totally convex, $G$-invariant subset of $M$, denoted by  $C_0$,  with nonempty boundary.  Clearly the boundary $\partial
C_0$ is also  $G$-invariant.  By inspecting the proof of Proposition 3.4 above we see that  the minimal nonempty superlevel set of the function $d_{\partial C_0}:=d(\cdot,
\partial C_0)$ on $C_0$ consists of a unique point, denoted by $S$.  $S$ is a soul of $M$. We claim that $S$ is fixed by the group $G$. The reason is very simple:  Otherwise there   exists $h \in G$ with $h\cdot S\neq S$. Then $d(S,\partial C_0)=d(S, h^{-1}\cdot \partial C_0)= d(h\cdot S, \partial C_0)< d(S,\partial C_0)$. A contradiction.  \hfill{$\Box$}

\vspace*{0.4cm}

The following proposition extends Theorems 3.9, 3.10 in
[CTZ].

\begin{prop} \label{prop 3.6} For every $\varepsilon>0$  there exist constants $C_1=C_1(\varepsilon)$ and $C_2=C_2(\varepsilon)$,
such that for every (restricted) ancient $\kappa$-orbifold solution
$(\mathcal{O},g(\cdot))$, for each space-time point $(x,t)$, there
is a radius $r$, $\frac{1}{C_1}(R(x,t))^{-\frac{1}{2}}< r<
C_1(R(x,t))^{-\frac{1}{2}}$, and an open neighborhood $B$,
$\overline{B(x,t,r)}\subset B \subset B(x,t, 2r)$, which falls into
one of the following categories:

(a) $B$ is an strong $\varepsilon$-neck centered at $(x,t)$,

(b) $B$ is an $\varepsilon$-cap,

(c) $B$ is diffeomorphic to a spherical orbifold
$\mathbb{S}^4//\Gamma$  (for a finite subgroup $\Gamma$ of $O(5)$).

Moreover, the scalar curvature in $B$  at time
$t$ is between $C_2^{-1}R(x,t)$ and $C_2R(x,t)$.
 \end{prop}

\noindent {\bf Proof} \ \ (Clearly we may assume $\varepsilon>0$ is sufficiently small.)  The arguments in [CTZ] (see proof of Theorems 3.9, 3.10 there) for the special case
when $\mathcal{O}$ has at most isolated singularities   do not work at certain steps in our more general situation.
So while our proof here is similar to that of Theorems 3.9, 3.10 in [CTZ] in many aspects, we have to modify it at some points. We only need to consider the case that the orbifold
$\mathcal{O}$ is noncompact and has positive curvature operator, since  the
other cases are clear (by Proposition 3.3 and the orbifold version of the main result of [H86] in [BW] ). Then by Proposition 3.4, $\mathcal{O}$ is
diffeomorphic to $\mathbb{R}^4// \Gamma$ for some finite subgroup
$\Gamma$ of $O(4)$.  We pull back the Ricci flow
$(\mathcal{O}, g(\cdot))$ to $(\mathbb{R}^4, \tilde{g}(\cdot))$, which is then a
$\Gamma$-invariant ancient $\kappa$-solution on a smooth manifold.  We fix  time $t=0$. By the proof of Proposition 3.4 in [CZ2], there is a point
$x_0\in \mathbb{R}^4$, such that for any given small $\eta>0$ there
is a positive constant $D(\eta)$, such that any point $x\in
\mathbb{R}^4$ with $R(x_0,0)d_{0}(x,x_0)^2\geq D(\eta)$ is the
center of a strong $\eta$-neck.  We rescale the solution $\tilde{g}(\cdot)$ so that
$R(x_0, 0)=1$ after the scaling; the rescaled solution will still be denoted by $\tilde{g}(\cdot)$. As in the proof of Theorem 3.9 in
[CTZ], we
use the manifold case of Proposition 2.3 (that is essentially Theorem C2.4 in [H97]) repeatedly to get Hamilton's canonical parametrization $\Phi:
\mathbb{S}^3 \times (A,B) \rightarrow  \mathbb{R}^4 $ such that the image of $\Phi$ contains $\mathbb{R}^4 \setminus
B(x_0,0, \sqrt{D(\eta)}+1)$, where $B(x_0, 0,\sqrt{D(\eta)}+1)$ denotes the  ball of radius $\sqrt{D(\eta)}+1$ w.r.t. the rescaled pull-back metric at time 0 centered at $x_0$.

Since  the ancient $\kappa$-solution  $(\mathbb{R}^4, \tilde{g}(\cdot))$  has asymptotic scalar curvature ratio $\infty$ by Perelman [P1], one can see that $(\mathbb{R}^4, \tilde{g}(0))$ ``splits off a line at infinity'' by using, for example, Proposition 2.3 in [CZ2]. In [CTZ] an argument using this property of  $(\mathbb{R}^4, \tilde{g}(0))$ shows that in effect,
the group $\hat{\Gamma}:=\Phi^{-1}\Gamma \Phi:=\{\Phi^{-1}\gamma \Phi| \hspace{1mm} \gamma \in \Gamma\}$  only acts on the $\mathbb{S}^3$ factor of $\mathbb{S}^3   \times (A,B)$, and the
  parametrization $\Phi$ can
be pushed down  to give Hamilton's canonical parametrization $\phi: \mathbb{S}^3//
\Gamma \times (A,B) \rightarrow \mathcal{O}$ of a part of
$\mathcal{O}$.

By the Fact showed just before this proposition, there exists a point in $\mathbb{R}^4$ which is a soul of $(\mathbb{R}^4, \tilde{g}(0))$ and which is fixed by $\Gamma$. We denote this point by $O$.

We claim that  $O$ has distance $\leq \sqrt{D(\eta)}+1$ from $x_0$
at time 0  when $\eta> 0$ is sufficiently small. The reason is as the following: Otherwise $O$ would be the
center of an $\eta$-neck.  Consider the distance function
$d_O(\cdot):=d(O, \cdot,0)=d_0(O, \cdot)$ from $O$ at time 0 (w.r.t. the rescaled pull-back metric). We know that
$d_O(\cdot)$ has no critical points in $\mathbb{R}^4 \setminus
\{O\}$, (since $O$ is a soul
of $\mathbb{R}^4$,) and the level sets $d_O^{-1}(s)$ are all homeomorphic to
$\mathbb{S}^3$  for all $s>0$, and in particular, are connected. But
on the other hand, if $O$ were the center of an $\eta$-neck,  the
level set $d_O^{-1}(s)$  would be disconnected for certain $s>0$
 (when $\eta>0$ is sufficiently small). A contradiction.

Let $S \in \mathcal{O}$ be the image of $O$. Now for any point $x \in \mathcal{O}$ with $d_0(x, S)\geq 2\sqrt{D(\eta)}+1$, and with $\eta> 0$ sufficiently small, we can use the $\phi$ above to give an $\varepsilon$-neck centered at $x$.

Then we want to find a cap in $\mathcal{O}$.  Let $\tilde{x} \in \mathbb{R}^4$ be a point with $d_0(O,
\tilde{x})=10\sqrt{D(\eta)}$ with $\eta> 0$ sufficiently small. Since $d_0(O, x_0)\leq \sqrt{D(\eta)}+1$, we see that $d_0(\tilde{x}, x_0)> 8\sqrt{D(\eta)}$. We denote the (``horizontal'') constant mean curvature
3-sphere passing through $\tilde{x}$ by $\Sigma$. $\Sigma$ is the center sphere of an $\varepsilon$-neck, so by
Theorem G1.1 in [H97], $\Sigma$ bounds an open subset $\Omega$ (in
$\mathbb{R}^4$) diffeomorphic to the standard open unit ball
$\mathbb{B}^4$ (hence also to $\mathbb{R}^4$) when $\varepsilon>0$ is sufficiently small.  From above we know that $\Omega$ is $\Gamma$-invariant, and $\Omega // \Gamma$ contains an $\varepsilon$-neck as its end . We want to modify the proof of Theorem G1.1 in [H97] to show that $\Omega$ is
$\Gamma$-equivariantly diffeomorphic to $\mathbb{R}^4$, thus $\Omega // \Gamma$ is diffeomorphic to the cap $\mathbb{R}^4 // \Gamma$.  (While at this point Chen-Tang-Zhu [CTZ] work  downstairs, it is more convenient for us to work upstairs equivariantly since our  necks may contain orbifold singularities. But the basic idea in both approaches is roughly the same as that in the proof of Theorem G1.1 in [H97].)

 Let $f=\pi_{(A,B)}\circ \Phi^{-1}:\mathbb{R}^4 \setminus B(x_0, 0,
\sqrt{D(\eta)}+1)   \rightarrow (A,B)$, where $\pi_{(A,B)}: \mathbb{S}^3
\times (A,B)  \rightarrow (A,B)$ is the natural projection. Consider   the  Busemann function $F$ on ($\mathbb{R}^4, \tilde{g}(0)$) constructed by using the base point  $O$, as in the proof of Proposition 3.4.   Clearly $F$ is $\Gamma$-invariant. From Appendix G in Hamilton [H97] we see  that  when $\eta >0$ is sufficiently small, in $\mathbb{R}^4 \setminus B(O, 0,
3\sqrt{D(\eta)})$
 the function $F$ is close to $f$ (if we multiply $f$ by $-1$ when necessary and shift it by a suitable constant) both in the supremum sense and in the Lipschitz sense. We smooth $F$ a little to $\tilde{F}$ by using, for example, the heat equation, such that $\tilde{F}$ is  strictly concave, $\Gamma$-invariant, and $\tilde{F}$ is close to $f$ in the $C^1$-norm in $\mathbb{R}^4 \setminus B(O,0,
3.5\sqrt{D(\eta)})$; cf. for example, Appendix F in [H97] and the proof of Proposition 3.7 in [CTZ]. The equivariant Morse theory gives that the (strict) superlevel set $\tilde{F}^{-1}(-c, +\infty)$ is $\Gamma$-equivariantly diffeomorphic to  $\mathbb{R}^4$ for large $c$, via the gradient flow of the function $\tilde{F}$.
 For any fixed $k\geq 4$ and $\delta>0$, let $\xi: \mathbb{R}^4 \rightarrow [0,1]$ be a smooth $\Gamma$-invariant function on $\mathbb{R}^4$  which  is 0 on the ball $B(O,0, k\sqrt{D(\eta)})$, and is 1 outside of the ball $B(O,0, (k+\delta)\sqrt{D(\eta)})$. We construct  the smooth  $\Gamma$-invariant function $(1-\xi)\tilde{F}+\xi f$ on $\mathbb{R}^4$, whose gradient flow  gives the $\Gamma$-equivariant diffeomorphism of $\tilde{F}^{-1}(-c, +\infty)$ with $(\xi f)^{-1}(-c', +\infty)$ for large $c$ and $c'$. (Compare [CTZ].) Finally we get the desired $\Gamma$-equivariant diffeomorphism of $\Omega$ with $\mathbb{R}^4$, since $\Omega$ is some (strict) superlevel set of $\xi f$.

The curvature estimate on $\Omega// \Gamma$ follows from Proposition 3.5
(i). \hfill{$\Box$}

\vspace*{0.4cm}

\noindent {\bf Remark }  It is not difficult to extend Naber's
classification of noncompact shrinking four-solitons with
nonnegative curvature operator ([N]) to the case of orbifolds. This
makes it possible to give  an alternative approach to Proposition
3.6 more along the original lines of [P1][P2].

\section{Metric surgery for Ricci flow  on 4-orbifolds}

\hspace *{0.4cm}

\noindent {\bf Definition  }(cf. [BBM]) \ \  Given an interval $I\subset
\mathbb{R}$, an evolving Riemannian orbifold is a pair
$(\mathcal{O}(t),g(t))$ ($t \in I$), where $\mathcal{O}(t)$ is a
(possibly empty or disconnected) orbifold  and $g(t)$ is a
Riemannian metric on $\mathcal{O}(t)$.  We say that it is
piecewise $C^1$-smooth if there exists a discrete subset $J$ of $ I$,
 such that the
following conditions are satisfied:

i. On each connected component of $I\setminus J$, $t \mapsto
\mathcal{O}(t)$ is constant, and $t \mapsto g(t)$ is $C^1$-smooth;

ii. For each $t_0\in J$, $\mathcal{O}(t_0)=\mathcal{O}(t)$ for any
$t< t_0$ sufficiently close to $t_0$, and $t\mapsto g(t)$ is left
continuous at $t_0$;

iii. For each $t_0 \in J\setminus$  $\{$sup$I\}$, $t\mapsto
(\mathcal{O}(t),g(t))$ has a right limit at $t_0$, denoted by
$(\mathcal{O}_+(t_0),g_+(t_0))$.

\hspace *{0.4cm}

As in [BBM], a time $t\in I$ is regular if $t$ has a neighborhood in
$I$ where $\mathcal{O}(\cdot)$ is constant and $g(\cdot)$ is
$C^1$-smooth. Otherwise it is singular.  We also denote by $f_{max}$
and $f_{min}$ the supremum and infimum of a function $f$,
respectively, as in [BBM].

\hspace *{0.4cm}

  \noindent {\bf Definition  }(Compare [BBM], [Hu1])\ \  A piecewise $C^1$-smooth
  evolving Riemannian 4-orbifold $\{(\mathcal{O}(t), g(t))\}_{t \in I }$ is a
  surgical solution of the Ricci flow if it has the following
  properties.

  i. The equation $\frac{\partial g}{\partial t}=-2Ric$ is satisfied
  at all regular times;

  ii.  For each singular time $t_0$ one has $(a_1+a_2)_{min}(g_+(t_0))\geq
  (a_1+a_2)_{min}(g(t_0))$, $(c_1+c_2)_{min}(g_+(t_0))\geq
  (c_1+c_2)_{min}(g(t_0))$, and $R_{min}(g_+(t_0))\geq
  R_{min}(g(t_0))$;

  iii. For each singular time $t_0$ there is a locally finite collection
  $\mathcal{S}$ of disjoint embedded $\mathbb{S}^3//\Gamma$'s in $\mathcal{O}(t_0)$
  (where $\Gamma$'s are finite subgroups of $O(4)$), and an orbifold $\mathcal{O}'$ such that

  (a) $\mathcal{O}'$ is obtained from $\mathcal{O}(t_0)\setminus \mathcal{S}$ by
  gluing back $\bar{\mathbb{B}}^4//\Gamma$'s;

 (b)$\mathcal{O}_+(t_0)$ is a union of some connected components of $\mathcal{O}'$ and
 $g_+(t_0)=g(t_0)$ on $\mathcal{O}_+(t_0)\cap \mathcal{O}(t_0)$;

(c) Each component of $\mathcal{O}'\setminus \mathcal{O}_+(t_0)$ is
diffeomorphic to a spherical orbifold, or a neck, or a cap, or  an
orbifold connected sum of at most two spherical orbifolds.

\hspace *{0.4cm}

Motivated  by the properties of 4-dimensional ancient
$\kappa$-orbifold solutions (see Propositions 3.5 and 3.6) and the standard
solution ([CZ2, Corollary A.2], which can be easily adapted to the
case of orbifold standard solution (which will be defined later) via
lifting), following [P2] (compare [BBM], [CaZ], [KL1] and [MT]), [CZ2]
 and [CTZ], we introduce the notion of canonical neighborhood.

\hspace *{0.4cm}

\noindent {\bf Definition}  Let $\varepsilon$ and $C$ be positive constants. A
point $(x,t)$ in a surgical solution to the Ricci flow is said to
have an $(\varepsilon,C)$-canonical neighborhood if there is   an open
neighborhood $U$, $\overline{B(x,t,\sigma)} \subset U\subset
B(x,t,2\sigma)$ with $C^{-1}R(x,t)^{-\frac{1}{2}}<\sigma
<CR(x,t)^{-\frac{1}{2}}$, which falls into one of the following four
types:

(a) $U$ is a strong $\varepsilon$-neck with center $(x,t)$,

(b) $U$ is an $\varepsilon$-cap with center $x$ for $g(t)$,

(c) at time $t$, $U$ is diffeomorphic to a closed spherical orbifold
$\mathbb{S}^4//\Gamma$,

\noindent and if moreover, the scalar curvature in $U$ at time $t$  is between $C^{-1}R(x,t)$ and $CR(x,t)$, and satisfies
the derivative estimates
\begin{equation*}
|\nabla R|< C R^{\frac{3}{2}} \hspace*{8mm} \mbox{and} \hspace*{8mm}
|\frac{\partial R}{\partial t}|< C R^2,
\end{equation*}
  and for case (c), the
curvature operator of $U$ is positive, and the infimal sectional
curvature of $U$ is greater than $C^{-1}R(x,t)$.

\hspace *{0.4cm}

\noindent {\bf Remark 1} Our definition of canonical neighborhood is slightly
different from that in [CTZ]. We include the derivative estimates
for the scalar curvature in case (c) also (while [CTZ] does not) for convenience. Note that by Proposition 3.5 these derivative estimates
hold uniformly for all restricted ancient $\kappa$-orbifold solutions. We also
impose a sectional curvature  condition in case (c). Note that by
using orbifold coverings and arguing as in the proof of [KL1, Lemma
59.7], it is easy to see that this condition is reasonable.

\hspace *{0.4cm}

\noindent {\bf Remark 2} Note that by Propositions 3.5, 3.6  and
[CZ2, Corollary A.2] (as adapted to the case of orbifold standard
solution), for every  $\varepsilon>0$, there exists a positive
constant $C(\varepsilon)$ such that each point in any ancient
$\kappa$-orbifold solution or the orbifold standard solution has an
$(\varepsilon,C(\varepsilon))$-canonical neighborhood, except that
for the orbifold standard solution, an $\varepsilon$-neck may not be
strong.

\hspace *{0.4cm}

We choose $\varepsilon_0>0$ such that $\varepsilon_0<10^{-4}$ and
such that when $\varepsilon \leq 2\varepsilon_0$, both the result of
Proposition 2.3 here and the result of Theorem G1.1 in [H97] (with $k=[\varepsilon^{-1}], L=\varepsilon^{-1}$) hold true. Let $\beta:=\beta(\varepsilon_0)$ be the
constant given by Lemma 2.5.  Define
$C_0:=\max\{100\varepsilon_0^{-1}, 2C(\beta\varepsilon_0/2)\}$,
where $C(\cdot)$ is given in the Remark 2 above. Fix $c_0>0$. Let
$\varrho_0, \Psi_0, L_0, P_0, S_0$ be the constants given in Lemma
3.2 by setting $c=c_0$ and $K=1$.

\hspace *{0.4cm}

Now we consider some a priori assumptions, which consist of the
pinching assumption and the canonical neighborhood assumption.

\hspace *{0.4cm}

\noindent {\bf Pinching assumption}: Let
   $\varrho_0$, $\Psi_0$, $L_0$, $P_0$, $S_0$  be positive constants as given above. A  surgical solution to the Ricci flow
  satisfies the pinching assumption (with pinching constants $\varrho_0,\Psi_0,L_0, P_0, S_0$) if there hold
\begin{equation}\begin{aligned}
& a_1+\varrho_0>0, \hspace{8mm} c_1+\varrho_0>0, \\
&  \max\{a_3,b_3,c_3\}\leq \Psi_0(a_1+\varrho_0), \hspace{8mm}
\max\{a_3,b_3,c_3\}\leq \Psi_0(c_1+\varrho_0), \\
&  \mbox{and}  \\
& \frac{b_3}{\sqrt{(a_1+\varrho_0)(c_1+\varrho_0)}}\leq 1+\frac{L_0
e^{P_0t}}{\max\{\ln \sqrt{(a_1+\varrho_0)(c_1+\varrho_0)},S_0\}}
\end{aligned}\end{equation}
at all points and  times.

\hspace *{0.4cm}

\noindent {\bf Canonical neighborhood assumption}:  Let $\varepsilon_0$ and
$C_0$ be given as above.
 Let $r: [0,+\infty)\rightarrow (0,+\infty)$ be a non-increasing function. An evolving Riemannian 4-orbifold $\{(\mathcal{O}(t), g(t))\}_{t \in I}$
 satisfies the canonical neighborhood assumption  $(CN)_r$ if  any  space-time point $(x,t)$ with  $R(x,t)\geq
r^{-2}(t)$ has an  $(\varepsilon_0,C_0)$-canonical neighborhood.

\hspace *{0.4cm}

Bounded curvature at bounded distance is one of the key ideas in
Perelman [P1], [P2]; compare [MT, Theorem 10.2], [BBB$^+$, Theorem
6.1.1] and [BBM, Theorem 6.4]. 4-dimensional versions have appeared
in [CZ2] and [Hu1]. The following is a extension of the version in
[Hu1, Theorem B.1].

\begin{prop} \label{prop 4.1}  For each $c, \varrho, \Psi, L, P, S, A, C >0$ and
each $\varepsilon \in (0,2\varepsilon_0]$, there exists $Q=Q(c,
\varrho, \Psi,L, P, S,A, \varepsilon, C)>0$ and $\Lambda=\Lambda(c,
\varrho, \Psi, L, P,S, A, \varepsilon, C)>0$ with the following
property. Let $I=[a, b]$ ($0 \leq a < b < \frac{1}{2c})$ and
$\{(\mathcal{O}(t), g(t))\}_{t\in I}$ be a surgical solution with
uniformly positive isotropic curvature ($a_1+a_2\geq c$,
$c_1+c_2\geq c$), with bounded curvature, and satisfying the
pinching condition  (3.1) (with constants $\varrho, \Psi, L, P, S$).
Let $(x_0, t_0)$ be a space-time point such that:

1. $R(x_0, t_0)\geq Q$;

2. For each point $y\in B(x_0, t_0, AR(x_0, t_0)^{-1/2})$, if $R(y,
t_0)\geq 4R(x_0, t_0)$, then $(y, t_0)$ has an $(\varepsilon,
C)$-canonical neighborhood.

\noindent Then for any $y\in B(x_0, t_0, AR(x_0, t_0)^{-1/2})$, we
have
\begin{equation*}
\frac{R(y, t_0)}{R(x_0, t_0)}\leq \Lambda.
\end{equation*}
\end{prop}

\noindent {\bf  Proof} \ \ We will adapt the proof of [BBB$^+$, Theorem 6.1.1] and
[BBM, Theorem 6.4] to our situation, incorporating an argument  in the proof of Proposition 4.2 in
[CTZ]. We argue by contradiction. Suppose the result is not true.
Then there exist constants  $c, \varrho, \Psi, L, P, S, A, C >0$ and
 $\varepsilon \in (0,2\varepsilon_0]$, sequences $Q_k \rightarrow +\infty$, $\Lambda_k \rightarrow +\infty$, and a sequences of pointed  surgical
solutions $(\mathcal{O}(t), g(t), (x_k, t_k))$ ($0 \leq a\leq t \leq
b < \frac{1}{2c})$) with uniformly positive isotropic curvature
($a_1+a_2\geq c$, $c_1+c_2\geq c$), with bounded curvature and
satisfying the pinching condition  (3.1) (with constants $\varrho,
\Psi, L, P, S$), such that:

1. $R(x_k, t_k)\geq Q_k$;

2. for each point $y\in B(x_k, t_k, AR(x_k, t_k)^{-1/2})$, if $R(y,
t_k)\geq 4R(x_k, t_k)$, then $(y, t_k)$ has an $(\varepsilon,
C)$-canonical neighborhood;

3. for each $k$, there exists $z_k\in B(x_k, t_k, AR(x_k, t_k)^{-1/2})$ with

\begin{equation*}
\frac{R(z_k, t_k)}{R(x_k, t_k)}> \Lambda_k.
\end{equation*}

\noindent For each $k$, consider the parabolic rescaling
\begin{equation*}
\bar{g}_k(\cdot):=R(x_k,t_k)g_k(t_k+\frac{\cdot}{R(x_k, t_k)}).
\end{equation*}
We will adopt the convention in [BBB$^+$] and [BBM] to put a bar on the points when the relevant geometric quantities are computed w.r.t. the metric $\bar{g}_k$.

Define
\begin{equation*}
\rho:=\sup \{s>0| \exists C(s)>0, \forall k\in \mathbb{N}, \forall\bar{y}\in B(\bar{x}_k,0,s), R(\bar{y},0)\leq C(s)\}.
\end{equation*}
It is easy to see that there exists, up to extracting a subsequence,  $\bar{y}_k \in B(\bar{x}_k,0,\rho)$ such that
\begin{equation*}
R(\bar{y}_k,0)\rightarrow +\infty \hspace*{8mm}  \mbox{and} \hspace*{8mm}   d_0(\bar{x}_k,\bar{y}_k)\rightarrow \rho \hspace*{8mm}  \mbox{as} \hspace*{8mm} k\rightarrow \infty.
\end{equation*}

We choose points $\bar{x}_k'$ and $\bar{y}_k'$ in  the geodesic segment $[\bar{x}_k\bar{y}_k]$ for large $k$ such
that $R(\bar{x}_k',0)=4C$, $R(\bar{y}_k',0)=R(\bar{y}_k,0)/(4C)$,
and $[\bar{x}_k'\bar{y}_k'] \subset [\bar{x}_k\bar{y}_k]$ is a
maximal subsegment on which
\begin{equation*}
4C\leq R(\cdot,0)\leq \frac{R(\bar{y}_k,0)}{4C},
\end{equation*}
with $\bar{x}_k'$ nearest to  $\bar{x}_k$.

As in [BBB$^+$] we can show that each point $\bar{z}$ in
$[\bar{x}_k'\bar{y}_k']$ has an $(\varepsilon, C)$ canonical
neighborhood which is a  strong $\varepsilon$-neck, say
$U(\bar{z})$, centered at $(\bar{z},0)$. Let $U_k$ be the union of
these $U(\bar{z})$'s. By using Proposition 2.3 repeatedly, the most part of $U_k$ (that is,  except for
the part near the two ends; one can give a more precise description of it with the aid of Lemma 2.2), denoted by $\mathcal{T}_k$, admits
Hamilton's canonical uniformization, say, $\Phi_k:\mathbb{S}^3\times
(A_k,B_k)\rightarrow \mathcal{T}_k$. Then
similarly as in the proof of [CTZ, Propositions 4.2, 4.4], we pull back
the rescaled solution $(\bar{g}_k(\cdot), \bar{x}_k')$ to $\mathbb{S}^3\times
(A_k,B_k)$ (which is a smooth manifold) via $\Phi_k$. The
pulled-back solutions (with the suitable base points)
sub-converge smoothly to a partial Ricci flow  (cf. [BBB$^+$]). Now
the rest of the argument is almost identical to that in [BBB$^+$,
Theorem 6.1.1]. For some of the details one can also consult Step 2
of proof of [CZ2, Theorem 4.1] (for the smooth (without surgery)
case) and Step 3 of proof of [CZ2, Proposition 5.4] (for the
surgical case).    \hfill{$\Box$}

\hspace *{0.4cm}

The following proposition extends [Hu1, Proposition 2.3]; compare
[BBM, Theorem 6.5] and [BBB$^+$, Theorem 6.2.1].

\begin{prop} \label{prop 4.2}\ \   Fix $c_0>0$.
 For any  $r$, $\delta>0$,
there exist    $h \in (0, \delta r)$ and $D> 10$, such that if
$(\mathcal{O}(\cdot),g(\cdot))$ is a complete surgical solution  with
uniformly positive isotropic curvature ($a_1+a_2\geq c_0$, $c_1+c_2
\geq c_0$), with bounded curvature, defined on a time interval
$[a,b]$ ($0\leq a < b < \frac{1}{2c_0}$) and satisfying the pinching
assumption and the canonical neighborhood assumption $(CN)_r$, then
the following holds:

 \noindent Let $t \in [a,b]$ and  $x, y, z \in \mathcal{O}(t)$ such that $R(x,t) \leq
2/r^2$, $R(y,t)=h^{-2}$  and $R(z,t)\geq D/h^2$. Assume there is a
curve $\gamma$ in $\mathcal{O}(t)$ connecting $x$ to $z$ and containing $y$, such that
each point of $\gamma$ with scalar curvature in $[2C_0r^{-2},
C_0^{-1}Dh^{-2}]$ is the center of an $\varepsilon_0$-neck.  Then
$(y,t)$ is the center of a strong $\delta$-neck.
\end{prop}

{\bf Proof}\ \ We essentially follow the proof of [BBM, Theorem 6.5]
and [BBB$^+$, Theorem 6.2.1] with some necessary modifications. (Compare [P2, Lemma 4.3], [CTZ,
Proposition 4.4].) We argue by contradiction. Otherwise, there exist
$r, \delta>0$,  sequences $h_k \rightarrow 0$, $D_k\rightarrow
+\infty$, a sequence of complete surgical solutions
$(\mathcal{O}_k(\cdot),g_k(\cdot))$  with bounded curvature and with uniformly
positive isotropic curvature ($a_1+a_2\geq c_0$, $a_1+a_2 \geq c_0$)
satisfying the  pinching assumption (with constants $\varrho_0,
\Psi_0, L_0, P_0, S_0$) and $(CN)_r$, and sequences $0 < t_k<
\frac{1}{2c_0}$, $x_k, y_k, z_k \in \mathcal{O}_k(t_k)$ with $R(x_k,t_k)\leq 2r^{-2}$, $R(y_k,t_k)=h_k^{-2}$
and $R(z_k,t_k)\geq D_kh_k^{-2}$, and finally a sequence of curves
$\gamma_k$ in $\mathcal{O}_k(t_k)$ connecting $x_k$ to $z_k$ and containing $y_k$, whose points
of scalar curvature in $[2C_0r^{-2},C_0^{-1}D_kh_k^{-2}]$ are
centers of $\varepsilon_0$-necks, but  $y_k$  is not  the center of a strong
$\delta$-neck.

Let $\bar{g}_k(t)=h_k^{-2}g_k(t_k+h_k^2t)$ for each $k$. For any $\rho>0$,  as in the proof of [BBM, Theorem 6.5], when $k$ is sufficiently large, $x_k \notin B(\bar{y}_k,\rho)$, $z_k \notin B(\bar{y}_k,\rho)$, and   $B(\bar{y}_k,\rho)$ is contained in the union of some $\varepsilon_0$-necks.  Now let $\Phi_k:\mathbb{S}^3\times (A_k,B_k)\rightarrow \mathcal{T}_k \subset \mathcal{O}_k(t_k)$ be Hamilton's canonical uniformization 
 whose image contains $y_k$ and which is maximal.
  Then
we pull back the parabolically rescaled solutions  $\bar{g}_k(t)$ to
$\mathbb{S}^3\times (A_k,B_k)$  via $\Phi_k$. The rest of the proof is almost the
same as in that of [BBM, Theorem 6.5] and [BBB$^+$, Theorem 6.2.1],
using Proposition 4.1 and local compactness theorem for Ricci flow. See also the proof of [Hu1, Proposition 2.3]. (By the way, in the definition of $\tau_0$ in the  proof of [Hu1, Proposition 2.3], $k(\rho)$ should read $k(\rho, \tau)$.)
\hfill{$\Box$}

\hspace *{0.4cm}

{\bf Remark}  For the purpose of this paper it suffices to use the following slightly weaker result which was actually proved in the previous versions of this paper:      Fix $c_0>0$.
 For any  $r$, $\delta>0$,
there exist    $h \in (0, \delta r)$ and $D> 10$, such that if
$(\mathcal{O}(\cdot),g(\cdot))$ is a complete surgical solution  with
uniformly positive isotropic curvature ($a_1+a_2\geq c_0$, $c_1+c_2
\geq c_0$), with bounded curvature, defined on a time interval
$[a,b]$ ($0\leq a < b < \frac{1}{2c_0}$) and satisfying the pinching
assumption and the canonical neighborhood assumption $(CN)_r$, then
the following holds:

 \noindent Let $t \in [a,b]$ and  $x, z \in \mathcal{O}(t)$ such that $R(x,t) \leq
2/r^2$  and $R(z,t)\geq D/h^2$. Assume there is a
curve $\gamma$ in $\mathcal{O}(t)$ connecting $x$ to $z$, such that
each point of $\gamma$ with scalar curvature in $[2C_0r^{-2},
C_0^{-1}Dh^{-2}]$ is the center of an $\varepsilon_0$-neck.  Then there is a point $y \in \gamma$ with $R(y,t)=h^{-2}$ such that
$(y,t)$ is the center of a strong $\delta$-neck.

\hspace *{0.4cm}

\noindent {\bf Remark}  It is clear from the above proof that the parametrization associated
to the  strong $\delta$-neck in the conclusion of Proposition 4.2 can be chosen to be Hamilton's canonical parametrization.

\hspace *{0.4cm}

Now we describe  Hamilton's metric surgery procedure [H97], which was
adapted to the case of orbifolds with isolated singularities in [CTZ]; we'll adapt it further to the case of more general orbifolds. We'll follow [CaZ], [CZ2],
and [CTZ]. First we describe the model surgery on the standard
cylinder, and define the orbifold standard solution. Consider the
semi-infinite cylinder $N_0=(\mathbb{S}^3//\Gamma) \times
(-\infty,4)$ with the standard metric $\bar{g}_0$ of scalar
curvature 1, where $\Gamma$ is a finite subgroup of isometries of
$\mathbb{S}^3$. Let $f$ be certain smooth nondecreasing convex function on $(-\infty, 4)$ as chosen in [CaZ, (7.3.7)] and [CZ2].
 Replace the
standard metric  $\bar{g}_0$ on the subspace $(\mathbb{S}^3//\Gamma)
\times [0,4)$ in $N_0$ by $e^{-2f}\bar{g}_0$.   The resulting metric
will induce  a complete metric (denoted by) $\hat{g}$ on the
 cap  $\mathbb{R}^4//\Gamma$.   We call an Ricci flow with initial data
$(\mathbb{R}^4//\Gamma,\hat{g})$ and with bounded curvature in any compact subinterval of $[0,
\frac{3}{2})$ the orbifold standard solution, which exists on the
time interval $[0, \frac{3}{2})$. Note that when $\Gamma=\{1\}$,
$(\mathbb{R}^4//\Gamma,\hat{g}(\cdot))$ is actually the  smooth standard
solution $(\mathbb{R}^4,\hat{g}(\cdot))$ with the initial metric
$(\mathbb{R}^4,\hat{g})$ constructed in [CZ2, Appendix A]. There is a
natural orbifold covering $\pi_\Gamma:
(\mathbb{R}^4,\hat{g})\rightarrow (\mathbb{R}^4//\Gamma,\hat{g})$, (we use the same $\hat{g}$). Denote by
$p_0$ the tip of the smooth standard solution, (which is the fixed
point of the $SO(4)$-action on
$(\mathbb{R}^4,\hat{g})$,)  and by $p_\Gamma=\pi_\Gamma(p_0) \in
\mathbb{R}^4//\Gamma$ the corresponding tip of the orbifold standard solution.

\noindent Then we describe a similar surgery procedure for the general case.
 Suppose we have a $\delta$-neck $N$ centered at $x_0$ in a Riemannian 4-orbifold $(\mathcal{O},g)$.
  Sometimes we will call $R^{-\frac{1}{2}}(x_0)$ the  radius of this neck. The surgery is to cut off the $\delta$-neck along
the center and glue back two orbifold caps $\mathbb{R}^4//\Gamma$
separately.  (Actually we first do surgery upstairs equivariantly,
then push it down.)
  Assume the center of the $\delta$-neck $N$ has $\mathbb{R}$
coordinate $z=0$.
 We construct a new metric on the glued back orbifold cap  (say on the left hand side) as
follows,
\begin{equation*}
 \tilde{g}= \begin{cases}
    g, & {z=0;}  \\
    e^{-2f}g, & {z\in [0,2];} \\
    \varphi e^{-2f}g+(1-\varphi)e^{-2f}h^2\bar{g}_0, & {z\in [2,3];} \\
    e^{-2f}h^2\bar{g}_0, & {z\in [3,4],}
  \end{cases}
\end{equation*}
where $\varphi$ is a smooth bump function with $\varphi=1$ for
$z\leq 2$, and $\varphi=0$ for $z\geq 3$, $h=R^{-\frac{1}{2}}(x_0)$,
and $\bar{g}_0$ is as above.
 We also use the same construction
on the right hand side with parameters $\bar{z}\in [0,4]$
($\bar{z}=8-z$).

 The following  lemma of Hamilton justifies the pinching assumption
of surgical solution.

 \begin{lem} \label{lem 4.3}\ \ (Hamilton [H97,D3.1]; compare [CZ2, Lemma 5.3], [CTZ, Lemma 4.3]) There
exist universal positive constants $\delta_0$, and a
constant $h_0$ depends only on $c_0$, such that given any surgical
solution with uniformly positive isotropic curvature ($a_1+a_2\geq
c_0$, $c_1+c_2\geq c_0$), satisfying the pinching assumption,
defined on $[a, t_0]$ ($0\leq a< t_0< \frac{1}{2c_0}$), if we
perform Hamilton's  surgery as described above at a $\delta$-neck
(if it exists) of radius $h$ at time $t_0$ with $\delta <\delta_0$
and $h \leq h_0$, then after the surgery, the pinching assumption
still holds at all points at time $t_0$, and any metric ball of radius $\delta^{-\frac{1}{2}}h$ with
center near the tip (i.e. the origin of the attached cap) is, after
scaling with the factor $h^{-2}$, $\delta^{\frac{1}{2}}$-close to
the corresponding ball of $(\mathbb{R}^4//\Gamma,\hat{g})$ for some
finite subgroup $\Gamma< O(4)$.

Furthermore, the  pinching assumption  holds after the surgery time
$t_0$ (by inspecting Hamilton's proof of pinching estimates in [H97,
Section B]).
\end{lem}

Usually we'll be given two non-increasing step functions $r,
\delta: [0,+\infty)\rightarrow (0, +\infty)$ as surgery parameters.
Let $h(r,\delta), D(r,\delta)$ be the associated parameters  as
determined in Proposition 4.2, which also depend on $\varepsilon_0$ and $C_0$, ($h$ is sometimes called the surgery
scale,) and let $ \Theta:=2Dh^{-2}$ be the curvature threshold for
the surgery process ( as in [BBM] and [Hu1]), that is, we'll do
surgery when $R_{max}(t)$ reaches $\Theta(t)$. Now we adapt two more
definitions from [BBM] and [Hu1].

\hspace *{0.4cm}

\noindent {\bf Definition} (compare [BBM], [Hu1] )\ \ Given an interval
$I\subset [0,+\infty)$, fix surgery parameters $r$, $\delta:
I\rightarrow (0,+\infty)$ (two non-increasing functions) and let
$h$, $D$, $\Theta=2Dh^{-2}$ be the associated cutoff parameters. Let
$(\mathcal{O}(t),g(t))$ ($t \in I$) be an evolving Riemannian
4-orbifold. Let $t_0 \in I$ and $(\mathcal{O}_+,g_+)$ be a
(possibly empty) Riemmanian 4-orbifold. We say that
$(\mathcal{O}_+,g_+)$ is obtained from
$(\mathcal{O}(\cdot),g(\cdot))$ by $(r,\delta)$-surgery at time
$t_0$ if

i. $R_{max}(g(t_0))=\Theta(t_0)$, and there is a locally finite
collection
  $\mathcal{S}$ of disjoint embedded $\mathbb{S}^3//\Gamma$'s in $\mathcal{O}(t_0)$ which are in the middle  of
  strong $\delta(t_0)$-necks with radius equal to the surgery scale $h(t_0)$, such that
  $\mathcal{O}_+$ is obtained from $\mathcal{O}(t_0)$ by doing
  Hamilton's surgery along these necks as described above  ( where $\Gamma$'s are  finite subgroups of $O(4)$
 ), and removing each of the following components:

  (a)  a component diffeomorphic to a spherical orbifold, with positive curvature operator and has sectional curvature bounded below by $C_0^{-1}R(x,t_0)$, where $x$ is some point in this component,

  (b)   a component diffeomorphic to a neck, and is  covered by (infinite many)
$\varepsilon_0$-necks,

 (c)  a component diffeomorphic to a cap,  and  is covered by an $\varepsilon_0$-cap and (infinite many) $\varepsilon_0$-necks,

 (d) a component diffeomorphic to an orbifold connected sum of at most two spherical orbifolds, and is covered by two $\varepsilon_0$-caps and /or (finite many)
$\varepsilon_0$-necks.

ii. If $\mathcal{O}_+\neq \emptyset$, then $R_{max}(g_+)\leq
\Theta(t_0)/2$.

\hspace *{0.4cm}

\noindent {\bf Definition} (cf. [BBM] and [Hu1])\ \  A surgical solution
$(\mathcal{O}(\cdot),g(\cdot))$ defined on some time interval
$I\subset [0,+\infty)$ is an $(r,\delta)$-surgical solution  if it has the  following
properties:

i.  It satisfies the pinching assumption, and $R(x,t) \leq \Theta (t)$ for all $(x,t)$;

ii. At each singular time $t_0\in I$,
$(\mathcal{O}_+(t_0),g_+(t_0))$ is obtained from
$(\mathcal{O}(\cdot),g(\cdot))$ by $(r,\delta)$-surgery at time
$t_0$;

iii. Condition $(CN)_r$ holds.

\hspace *{0.4cm}

 Recall that in our 4-dimensional case, $g(\cdot)$ is
$\kappa$-noncollapsed (for some $\kappa>0$) on the scale $r$ at time
$t$ if at any point $x$, whenever $|Rm|\leq r^{-2} \hspace{2mm} \mbox{on}
\hspace{2mm} P(x, t, r, -r^2) \hspace{2mm}$ we have $ \hspace{2mm}
\mbox{vol} B(x,t, r)\geq \kappa r^4$.
 Let $\kappa: I \rightarrow (0, +\infty)$ be a function.
We say $\{(\mathcal{O}(t), g(t))\}_{t \in I}$ has property $(NC)_\kappa$ if it
is $\kappa (t)$-noncollapsed  on all scales $\leq 1$ at any time
$t\in I$. An $(r,\delta)$-surgical solution which also satisfies
condition $(NC)_\kappa$  is called an $(r,\delta,\kappa)$-surgical
solution.

\hspace *{0.4cm}

The following proposition extends [Hu1, Proposition 2.7], and is
analogous to [BBM, Proposition A].

\begin{prop} \label{prop 4.4}\ \ Fix $c_0>0$. There exists a positive
constant $\tilde{\delta}$ (depending only on $c_0>0$) with the
following property: Let $r, \delta$ be surgery parameters, let
$\{(\mathcal{O}(t), g(t))\}_{t\in (a,b]}$ ( $0<a<b<\frac{1}{2c_0}$)
be an $(r, \delta)$-surgical solution  with uniformly positive
isotropic curvature ($a_1+a_2\geq c_0$, $c_1+c_2 \geq c_0$). Suppose
that $\delta\leq \tilde{\delta}$, and $R_{max}(b)=\Theta=\Theta(b)$.
 Then there exists a Riemannian orbifold
$(\mathcal{O}_+,g_+)$ which is obtained from $(\mathcal{O}(\cdot),g(\cdot))$ by
$(r,\delta)$-surgery at time $b$, such that

i.  $g_+$ satisfies the pinching assumption at time $b$;

ii. $(a_1+a_2)_{min}(g_+(b))\geq
  (a_1+a_2)_{min}(g(b))$, $(c_1+c_2)_{min}(g_+(b))\geq
  (c_1+c_2)_{min}(g(b))$, and $R_{min}(g_+(b))\geq R_{min}(g(b))$.
\end{prop}

\noindent {\bf Proof}\ \  Let  $\delta_0$ and $h_0$ be as  given in Lemma
4.3. Set $\tilde{\delta}=\frac{1}{2} \min \{c_0^{\frac{1}{2}}h_0,
\delta_0\}$. The idea is to consider a maximal collection $\{N_i\}$
of pairwise disjoint cutoff necks in $\mathcal{O}(b)$, whose
existence is guaranteed by Zorn's Lemma. (Here, following [BBM], a strong $\delta$-neck centered at some point
$(x,t)$ of scalar curvature $h^{-2}(t)$ is called a cutoff neck.) We want to show that such
a collection is locally finite.  Note that in the orbifold case we
need a new argument to guarantee this; compare the volume argument
in the manifold case in [BBM] and [Hu1]. We argue by contradiction.
Otherwise there is
 a sequence of cutoff necks (still denoted by $\{N_i\}$ ) with center
$y_i$ (with $R(y_i,b)=h^{-2}(b))$, where $N_i$ is diffeomorphic to $(\mathbb{S}^3//{\Gamma_i})
  \times \mathbb{I}$ with $|\Gamma_i|\rightarrow \infty$ as $i\rightarrow \infty$,
 and all $N_i$'s are contained in a compact subset $K$ of $\mathcal{O}(b)$.
Then there is a subsequence of $\{y_i\}$ (still denoted by $\{y_i\}$) which converges to a point $y$ in $K$.
  We have $R(y,b)=h^{-2}(b)$, so $(y,b)$ has a canonical neighborhood $U$,
 which is impossible, as can be seen as follows.

 i. If $U$ is in case (a) in the
definition of canonical neighborhood, we get a contradiction by
using the assumptions $|\Gamma_i|\rightarrow \infty$ and
$y_i\rightarrow y$ as $i\rightarrow \infty$ and Proposition 2.3.

ii. If $U$ is in case (b), we may assume that $(y,b)$ is not the
center of an $\varepsilon_0$-neck, otherwise we can argue as in case
i above and get a contradiction. Then, for $y_i$ sufficiently close
to $y$, $y_i$ can not be the center of  any strong $\delta$-neck.
Again a contradiction.

iii. If $U$ is in case (c), we get a contradiction by comparing
sectional curvature.

 Now with the help of Proposition 4.2, Lemma 4.3 and
Proposition 2.4, the rest of the argument is similar to that in the proof of [Hu1, Proposition
2.7], and is omitted.

\hfill{$\Box$}

\section{Existence of $(r, \delta, \kappa)$-surgical solutions}

As in [BBM], if $(\mathcal{O}(\cdot), g(\cdot))$ is a piecewise $C^1$ evolving
orbifold defined on some interval $I\subset \mathbf{R}$ and
$[a,b]\subset I$, the restriction of $g$ to $[a,b]$, still denoted
by $g(\cdot)$, is the evolving orbifold
\begin{equation*}
 t\mapsto \begin{cases}
     (\mathcal{O}_+(a), g_+(a)), &       t=a, \\
    (\mathcal{O}(t), g(t)), &          t\in (a,b]. \\
 \end{cases}
 \end{equation*}

The following proposition extends [Hu1, Proposition 3.1]. Compare
[P2, Lemma 4.5], [BBB$^+$, Theorem 8.1.2], [BBM, Theorem 8.1], [CaZ,
Lemma 7.3.6], [KL1, Lemma 74.1], [KL2, Lemma 7.29], [MT, Proposition 16.5] and [Z, Lemma
9.1.1], see also the formulation in the proof of [CZ2, Lemma 5.5].

\begin{prop} \label{prop 5.1}\ \  Fix $c_0>0$. For all $A>0, \theta \in
(0,\frac{3}{2})$ and $\hat{r}>0$, there exists
$\hat{\delta}=\hat{\delta}(A,\theta,\hat{r})>0$ with the following
property. Let $r(\cdot)\geq \hat{r}$, $\delta(\cdot) \leq
\hat{\delta}$ be two positive, non-increasing  step functions on
$[a,b)$ ($0\leq a< b<\frac{1}{2c_0}$), and let
$(\mathcal{O}(\cdot),g(\cdot))$ be a surgical solution  with
uniformly positive isotropic curvature ($a_1+a_2\geq c_0$, $c_1+c_2
\geq c_0$), defined on $[a,b]$, such that it satisfies the pinching
assumption on $[a,b]$, that $R(x,t)\leq \Theta(r(t), \delta(t))$ for
all space-time points with $t \in [a,b)$, that at any singular time
$t_0\in [a,b)$, $(\mathcal{O}_+(t_0), g_+(t_0))$ is obtained from
$(\mathcal{O}(\cdot), g(\cdot))$ by $(r, \delta)$-surgery, and that
any point $(x,t)$ ($t\in [a,b)$) with $R(x,t) \geq
(\frac{r(t)}{2})^{-2}$ has a $(2\varepsilon_0, 2C_0)$-canonical
neighborhood. Let $t_0\in [a,b)$ be a singular time. Consider the
restriction of $(\mathcal{O}(\cdot),g(\cdot))$ to $[t_0,b]$. Let
$p\in \mathcal{O}_+(t_0)$ be the tip of some surgery cap of scale
$h(t_0)$, and let $t_1=\min \{b,t_0+\theta h^2(t_0)\}$. Then either

(i) The parabolic neighborhood $P(p,t_0,Ah(t_0),t_1-t_0)$ is
unscathed, and is, after scaling with factor $h^{-2}(t_0)$ and
shifting time $t_0$ to zero, $A^{-1}$-close to
$P(p_\Gamma,0,A,(t_1-t_0)h^{-2}(t_0))$ (where $p_\Gamma$ is the tip
of the cap of the orbifold standard solution $\mathbb{R}^4//\Gamma$ for some finite subgroup
$\Gamma< O(4)$), or

(ii) The assertion (i) holds with $t_1$ replaced by some $t^+ \in
[t_0, t_1)$; moreover  $B(p,t_0,Ah(t_0))$ is removed by the surgery
at time $t^+$.
\end{prop}

We will follow the proof of [BBB$^+$, Theorem 8.1.2] and [BBM,
Theorem 8.1].

 Let $\mathcal{M}_0=(\mathbb{R}^4,\hat{g}(\cdot))$ be
the smooth standard solution, and  $0<T_0<\frac{3}{2}$.

\begin{prop} \label{prop 5.2}\ \ (Compare [BBB$^+$, Theorem 8.1.3] and [Hu1, Lemma
3.2]) For all $A, \Lambda>0$, there exists
$\rho=\rho(\mathcal{M}_0,A,\Lambda)>A$ with the following
properties. For a finite subgroup $\Gamma< O(4)$,  let $U$ be an open subset of $\mathbb{R}^4//\Gamma$,
and $T\in (0,T_0]$. Let $g(\cdot)$ be a Ricci flow defined on
$U\times [0,T]$, such that the ball $B(p_\Gamma,0,\rho)\subset U$ is
relatively compact. Assume that

(i) $||Rm(g(\cdot))||_{0,U\times [0,T],g(\cdot)}\leq \Lambda$,

(ii) $g(0)$ is $\rho^{-1}$-close to $\hat{g}(0)$ on
$B(p_\Gamma,0,\rho)$.

\noindent Then $g(\cdot)$ is $A^{-1}$-close to $\hat{g}(\cdot)$ on
$B(p_\Gamma,0,A)\times [0,T]$.

\hspace *{0.4cm}

\noindent Here, $||Rm(g(\cdot))||_{0,U\times [0,T],g(\cdot)}:=
\sup_{U \times [0,T]} \{|Rm_{g(t)}(x)|_{g(t)}\}$.
\end{prop}

\noindent {\bf Proof}\ \ We argue by contradiction. Otherwise there exist
$A,\Lambda >0$, and a sequence of Ricci flows $g_k(\cdot)$ defined
on $U_k \times [0,T_k]$ (where $U_k\subset \mathbb{R}^4//\Gamma_k$
with finite subgroup $\Gamma_k<O(4)$, and $T_k\leq T_0$), a sequence
$\rho_k\rightarrow +\infty$ as $k\rightarrow \infty$, such that
$B(p_k,0,\rho_k)\subset U_k$ are relatively compact, where
$p_k:=p_{\Gamma_k}$, and

(i) $|Rm_{g_k(t)}|_{g_k(t)}\leq \Lambda$ on $U_k \times [0,T_k]$,

(ii) $g_k(0)$ is $\rho_k^{-1}$-close to $\hat{g}(0)$ on
$B(p_k,0,\rho_k)$,

\noindent but for some $t_k \in [0,T_k]$, $g_k(t_k)$ is not
$A^{-1}$-closed to $\hat{g}(t_k)$ on $B(p_k,0,A)$.

We pull back the solutions  $g_k(\cdot)$ (and $\hat{g}(\cdot)$) to
$\mathbb{R}^4$ via $\pi_k:=\pi_{\Gamma_k}: \mathbb{R}^4 \rightarrow
\mathbb{R}^4//\Gamma_k$. Then we may acquire that

$|Rm_{\pi_k^*g_k(t)}|_{\pi_k^*g_k(t)}\leq \Lambda$ on
$\pi_k^{-1}(U_k) \times [0,T_k]$,

$\pi_k^*g_k(0)$ is $\rho_k^{-1}$-close to $\hat{g}(0)$ on
$B(p_0,0,\rho_k)$, but

$\pi_k^*g_k(t_k)$ is not $A^{-1}$-close to $\hat{g}(t_k)$ on
$B(p_0,0,A)$.

 \noindent Now we can argue as in [BBB$^+$], using a stronger  version of Shi's derivative
 estimates  ([LT,Theorem 11], see also
[MT, Theorem 3.29]),
  Hamilton's compactness theorem for Ricci flow ([H95]) and Chen-Zhu's
  uniqueness theorem  for complete Ricci flow ([CZ1]), to get a contradiction.
\hfill{$\Box$}

\begin{cor} \label{cor 5.3}\ \ (Compare [BBM, Corollary 8.3]  and [Hu1,
Corollary 3.3]) Let $A>0$. There exists $\rho=\rho(\mathcal{M}_0,A)
>A$ with the following properties. Let $\{(\mathcal{O}(t),g(t))\}_{t\in
[0,T]}$ ($T \leq T_0$) be a surgical solution. Assume that

(i) $(\mathcal{O}(\cdot),g(\cdot))$ is a parabolic rescaling  of
some surgical solution which satisfies the pinching assumption,

(ii)  $|\frac{\partial R}{\partial t}|\leq 2C_0 R^2$ at any
space-time point $(x,t)$ with $R(x,t)\geq 1$.

\noindent Let $p\in \mathcal{O}_+(0)$ and $t_0\in (0,T]$ be such that

(iii) $B(p,0,\rho)$ is $\rho^{-1}$-close to $B(p_\Gamma,0,\rho)
\subset \mathbb{R}^4//\Gamma$  for some finite subgroup $\Gamma$ of $O(4)$,

(iv) $P(p,0,\rho,t_0)$ is unscathed.

\noindent Then $P(p,0,A,t_0)$ is $A^{-1}$-close to
$P(p_\Gamma,0,A,t_0)$.
\end{cor}

\noindent {\bf Proof} The proof is similar to that of Corollaries 8.2.2 and
8.2.4 in [BBB$^+$].  \hfill{$\Box$}

\hspace *{0.4cm}

 Then to finish the proof of Proposition 5.1, we can proceed as in the proof of  [BBB$^+$, Theorem 8.1.2] and [BBM, Theorem 8.1], using
Corollary 5.3.

\hspace *{0.4cm}

The following lemma extends [Hu1, Lemma 3.5], and guarantees  the
non-collapsing under a slightly weaker canonical neighborhood
assumption. Compare [P2, Lemma 5.2],  [CTZ, Lemma 4.5], [ KL1, Lemma
79.12] and [BBM, Proposition C].

\begin{lem} \label{lem 5.4}\ \ Fix $c_0>0$. Suppose $0< r_- \leq \varepsilon_0
$, $\kappa_->0$, and $0< E_-<
 E< \frac{1}{2c_0}$. Then there exists
 $\kappa_+=\kappa_+(r_-,\kappa_-,E_-,E)>0$, such that for any  $r_+$, $0< r_+ \leq r_-$, one can find
  $\delta'=\delta'( r_-,r_+,\kappa_-,E_-,E)>0 $,  with the following
 property.

 \noindent Suppose that $0\leq a<b<d < \frac{1}{2c}$, $b-a\geq E_-$, $d-a\leq E$. Let $r$ and $\delta
$ be two positive, non-increasing  step functions on $[a,d)$ with
$\varepsilon_0\geq r\geq r_-$ on $[a,b)$,
 $\varepsilon_0\geq r\geq r_+$ on $[b,d)$  and $\delta \leq \delta'$ on $[a,d)$.
Let  $(\mathcal{O}(\cdot),g(\cdot))$ be a surgical solution
 with
uniformly positive isotropic curvature ($a_1+a_2\geq c_0$, $c_1+c_2
\geq c_0$), defined on the time interval $[a, d]$, such that it
satisfies the pinching assumption on $[a, d]$, that $R(x,t)\leq
\Theta(r(t),\delta(t))$ for all space-time points $(x,t)$ with $t\in [a,d)$,
that at any singular time $\tilde{t}\in [a,d)$, $(\mathcal{O}_+(\tilde{t}),
g_+(\tilde{t}))$ is obtained from $(\mathcal{O}(\cdot), g(\cdot))$ by $(r,
\delta)$-surgery, that the conditions $(CN)_{r}$ and
$(NC)_{\kappa_-}$ hold on $[a,b)$, and that any point $(x,t)$ ($t\in
[b, d)$) with $R(x,t) \geq (\frac{r(t)}{2})^{-2}$ has a
$(2\varepsilon_0, 2C_0)$-canonical neighborhood. Then
$(\mathcal{O}(\cdot),g(\cdot))$ satisfies $(NC)_{\kappa_+}$ on
$[b,d]$.
\end{lem}

\noindent {\bf Proof} \ \ W.l.o.g. we may assume $r(\cdot)=r$ is constant on $[b,d)$.
Note that the property $(NC)_{\kappa_+}$ is closed w.r.t. time. Fix any $0< r_0 \leq 1$. Let
$t_0\in [b,d)$ and $x_0\in \mathcal{O}(t_0)$. Assume
$|Rm(\cdot,\cdot)| \leq r_0^{-2}$ on $P(x_0, t_0,r_0, -r_0^2)$, we
want to bound $vol_{t_0}(B(x_0,t_0,r_0))/r_0^4$ from below. We
will consider three cases.

Case 1. $r_0\geq \frac{r}{C(\varepsilon_0)}$, where $C(\varepsilon_0)$ is a constant depending only on $\varepsilon_0$, and is to be determined in Case 2. Using Perelman's reduced volume [P1, P2] and Proposition 5.1,
 the argument of Step 1 in the proof of [CTZ, Lemma 4.5] can be adapted to our case
without essential changes.

Case 2.  $r_0 < \frac{r}{C(\varepsilon_0)}$  and  there is a point $x$ in
the connected component containing $x_0$ such that $R(x,t_0)<
4r^{-2}$. By an argument using orbifold Bishop-Gromov theorem ([B93][Lu]) we may assume w.l.o.g. that there is some point $(x', t') \in \overline{P(x_0, t_0,r_0, -r_0^2)}$ such that $|Rm(x', t')|=r_0^{-2}$. (Compare Lemma 10.1.2 in  [BBB$^+$], and p.232-233 in [CZ2].)  Then we can get $R(x_0, t_0)\geq 4r^{-2}$ by choosing $C(\varepsilon_0)$ sufficiently large. As in Step 2 in the proof of [CTZ, Lemma 4.5], using the assumption that there is a point $x$ in
the connected component containing $x_0$ such that $R(x,t_0)<
4r^{-2}$, one can utilize  Hamilton's
canonical parametrization to find a tube, such that one end of the tube is adjacent to a neck or cap centered at $(x_0, t_0)$, and the scalar curvature $\leq C(\varepsilon_0)^2r^{-2}$ near the other end. Then we can reduce this case to Case 1.

Case 3.   $r_0 < \frac{r}{C(\varepsilon_0)}$  and  every point $x$ in the
connected component containing $x_0$ has $R(x,t_0)\geq 4r^{-2}$.
Then by our assumption at time $t_0$ this component, denoted by $K$,
is covered by $(2\varepsilon_0, 2C_0)$-canonical neighborhoods,  and
its topology can be described with the help of  Proposition 2.4.  By  definition of canonical neighborhood, to  bound $vol_{t_0}(B(x_0,t_0,r_0))/r_0^4$ from below we only need   to
uniformly control the topology of this component  $K$.  We consider
further two subcases.

Subcase a).  $K \times [b,t_0]$ is unscathed.  Note that
 $(x_0,t_0)$ has a $(2\varepsilon_0, 2C_0)$-canonical neighborhood, denoted by $U$. It suffices to uniformly control the topology of $U$.  We may assume
$\inf_{t \in [b,t_0]} R(x_0,t)< 4r^{-2}$, since otherwise the topology of $U$ is uniformly controlled by using the $\kappa_-$-noncollapsing at  time $b$. Now
let $t_0''= \sup \{t|~ b\leq t \leq t_0, R(x_0, t)=4r^{-2}\}$. By the time derivative estimate
in the definition of canonical neighborhood and pinching assumption,  we  can find $t_0' \in [t_0'', t_0]$ such that  $|Rm| \leq
C(\varepsilon_0)^2r^{-2}$ on $P(x_0, t_0', r, -r^2)$, and we are in a similar situation
as in Case 1, and we are done.  (Compare the argument for Case 2 in [BBM, Section 10.4].)

Subcase b).  $K \times [b,t_0]$ is scathed.  Let $t_1 \in [b, t_0)$
be the last singular time (w.r.t. $K$), then there is a
corresponding strong $\delta$-neck in $(\mathcal{O}(t_1), g(t_1))$,
 and   we only need to uniformly control the topology of this $\delta$-neck. If there is a point $x$ in the component $K_1$ of $\mathcal{O}(t_1)$ containing this $\delta$-neck
  with $R(x, t_1)<4r^{-2}$, then by Case 2 we are done.  Otherwise at time $t_1$ this component $K_1$ is covered by $(2\varepsilon_0, 2C_0)$-canonical neighborhoods.
  Now let $t_2 \in [b, t_1)$  be the last singular time (w.r.t. $K_1$), and we repeat the above argument.  After finite steps we will be in a similar situation as in Subcase a), and we are done.
\hfill{$\Box$}

\hspace *{0.4cm}

To justify the canonical neighborhood assumption needed,  we
 extend [Hu1, Proposition 3.6] to our more general situation here. The  argument is similar to that in the proof of [Hu1, Proposition 3.6] and  [CZ2, Proposition 5.4], with the aid of  a weak openness
(w.r.t. time) property of the canonical neighborhood condition in
the noncompact orbifold case (extending the noncompact manifold case
in [Hu1]), Lemma 2.5, Propositions 3.5, 3.6, 4.1, 5.1, Lemma
5.4, Hamilton's Harnack estimate [H93], the compactness theorem
for Ricci flow ([H95], [Lu], [KL1], [KL2], [To]), and Proposition 6.1 in the appendix of [CTZ].

\hspace *{0.4cm}

Once the canonical neighborhood assumption is justified, one can
prove the following theorem (extending [Hu1, Theorem 3.4]) similarly
as  in [Hu1] with the help of Lemma 3.2, Lemma 4.3, Proposition 4.4
and Lemma 5.4.

\begin{thm} \label{thm 5.5}\ \ Given  $c_0$, $v_0>0$, there are surgery parameter
sequences
\begin{equation*}
\mathbf{K}=\{\kappa_i\}_{i=1}^\infty, \hspace{2mm}
\Delta=\{\delta_i\}_{i=1}^\infty, \hspace{2mm}
\mathbf{r}=\{r_i\}_{i=1}^\infty
\end{equation*}
such that the following holds. Let $r(t)=r_i$ and
$\delta(t)=\delta_i$  and $\kappa(t)=\kappa_i$ on $[(i-1)2^{-5},
i\cdot2^{-5})$, $i=1, 2, \cdot\cdot\cdot$.  Let $(\mathcal{O},g_0)$
be  a complete 4-orbifold with uniformly positive isotropic
curvature ($a_1+a_2\geq c_0$, $c_1+c_2\geq c_0$), with $|Rm|\leq 1$,
and with vol $B(x,1)\geq v_0$ at any point $x$. Then there exists an
$(r,\delta, \kappa)$-surgical solution to the Ricci flow with initial data
$(\mathcal{O},g_0)$, which becomes extinct before the time
$\frac{1}{2c_0}$.
\end{thm}

\section {Proof of Theorems 1.1 and 1.3}

\begin{prop} \label{prop 6.1}   Let $\mathcal{F}$ be a class of closed
4-orbifolds. Let $\mathcal{O}$ be a 4-orbifold. Suppose there exists
a finite sequence of 4-orbifolds $\mathcal{O}_0$, $\mathcal{O}_1$,
$\cdot\cdot\cdot$, $\mathcal{O}_k$ such that
$\mathcal{O}_0=\mathcal{O}$, $\mathcal{O}_k=\emptyset$, and for each
$i$ ($1\leq i \leq k$), $\mathcal{O}_i$ is obtained from
$\mathcal{O}_{i-1}$ by cutting off along a locally finite collection
of pairwise disjoint, embedded spherical 3-orbifolds
$\mathbb{S}^3//\Gamma$'s, gluing back $\bar{\mathbb{B}}^4//\Gamma$'s,
and removing some components that are orbifold connected sums of
members of $\mathcal{F}$. Then each component of $\mathcal{O}$ is an
orbifold connected sum of members of $\mathcal{F}$.
\end{prop}

\noindent {\bf Proof}\ \  The proof is elementary, and is almost identical to
that of [BBM, Proposition 2.6], so we will omit it. \hfill{$\Box$}

\hspace *{0.4cm}

\noindent {\bf Proof of Theorem 1.3}.  Let $(\mathcal{O}, g_0)$ be a complete, connected
Riemannian 4-orbifold with uniformly positive isotropic curvature
and with bounded geometry. After normalization   $(\mathcal{O}, g_0)$ satisfies the condition in Theorem 5.5. By Theorem 5.5 we can construct
an $(r,\delta, \kappa)$-surgical solution to the Ricci flow starting with $(\mathcal{O}, g_0)$  which becomes extinct in finite time. Recall that each point in any component that is removed in the process of surgery is contained in a canonical neighborhood, so  any such component  is either diffeomorphic to a spherical 4-orbifold, or  covered by $\varepsilon_0$-caps and/ or $\varepsilon_0$-necks. In the latter case  such component  must appear in the list in the conclusion  of Proposition 2.4.  In any case, it is a
(possibly infinite) orbifold connected sum of spherical 4-orbifolds (by using an observation in Section 2).
Then by our surgery procedure and  Proposition 6.1, $\mathcal{O}$ is diffeomorphic to a (possibly
infinite) orbifold connected sum of spherical 4-orbifolds.  Now we
argue that the diffeomorphism types of the spherical 4-orbifolds that
appear in the orbifold connected sum decomposition of $\mathcal{O}$ is finite,
which will finish the proof of Theorem 1.3. We divide the analysis
into two cases.

i. Consider those components which  are removed in our surgery
procedure and each of which contains at least an
$\varepsilon_0$-neck. Our assumption on uniformly positive isotropic
curvature imply that the (``horizontal'') $\mathbb{S}^3//\Gamma$ cross section of
these $\varepsilon_0$-necks must have Ricci curvature uniformly
bounded below away from zero, which, combined with the
non-collapsing property and the boundedness of the sectional
curvature,
 implies that the isomorphism classes of $\Gamma$ are finite.

ii.  By definition of our surgery procedure, any other  component (that is, component which does not contain any $\varepsilon_0$-neck) which  is removed in our surgery
procedure must be diffeomorphic to a spherical orbifold
 (say $\mathbb{S}^4//\Gamma$, where
$\Gamma$ is a finite subgroup of $O(5)$), and  have  positive curvature operator with sectional
curvature bounded from below  by some positive constant. By the
orbifold Myers theorem (cf. for example [Lu]) the diameters of these
components are uniformly bounded above, which combined with the
non-collapsing property and the boundedness of sectional curvature
gives the desired finiteness of isomorphism classes of $\Gamma$ (cf.
[B]).

 \hfill{$\Box$}

\hspace *{0.4cm}

The following proposition deals with a special case of Proposition
2.4 more explicitly.

\begin{prop} \label{prop 6.2}\ \ Let $\varepsilon\in (0, 2\varepsilon_0]$.
Let $(\mathcal{O},g)$ be a complete, connected 4-orbifold with at most
isolated singularities. If each point of $\mathcal{O}$ is the center of an
$\varepsilon$-neck or an $\varepsilon$-cap, then $\mathcal{O}$ is
diffeomorphic to  a mapping torus $\mathbb{S}^3 / \Gamma {\times}_f \mathbb{S}^1$ (where $\Gamma$ is a  finite  subgroup
of isometries of $\mathbb{S}^3$ which acts freely on $\mathbb{S}^3$),
$\mathbb{S}^3 / \Gamma \times \mathbb{R}$, a smooth cap ($C_\Gamma^\sigma$ or $\mathbb{R}^4$), an orbifold cap of type I ($C_\Gamma$),
an orbifold cap of type II ($\mathbb{S}^4//(x,\pm
  x')\setminus \bar{\mathbb{B}}^4$), $C_\Gamma^\sigma {\cup}_f C_{\Gamma}^{{\sigma}'}$, $C_\Gamma^\sigma {\cup}_f C_{\Gamma}$,
 $C_\Gamma  {\cup}_f C_{\Gamma}$, $\mathbb{S}^4//(x,\pm x')$, $\mathbb{S}^4//(x,\pm x') \sharp
 \mathbb{R}\mathbb{P}^4$,  or $\mathbb{S}^4//(x,\pm x') \sharp \mathbb{S}^4//(x,\pm x')$.
\end{prop}

\noindent {\bf Proof}. The result in the compact case has been given in [CTZ] (see p. 61 and  p. 72 there). In general, one can argue as in [Hu1, Proposition
2.6], using Proposition 2.3.  Note that each of $\mathbb{S}^3 / \Gamma {\times}_f \mathbb{S}^1$, $C_\Gamma^\sigma {\cup}_f C_{\Gamma}^{{\sigma}'}$, $C_\Gamma^\sigma {\cup}_f C_{\Gamma}$ and
 $C_\Gamma  {\cup}_f C_{\Gamma}$ can be written as an orbifold connected sum of at most two spherical orbifolds (with at most isolated singularities), which was already observed in [CTZ], see also the proof of Proposition 2.4.
\hfill{$\Box$}

\hspace *{0.4cm}

\noindent {\bf Proof of Theorem 1.1}. Let $(X, g_0)$ satisfy the assumption of
Theorem 1.1. After normalization   $(X, g_0)$ satisfies the condition in Theorem 5.5. By Theorem 5.5 we can construct
an $(r,\delta, \kappa)$-surgical solution to the Ricci flow starting with $(X, g_0)$  which becomes extinct in finite time. Any component that is removed in the process of surgery is an orbifold with at most isolated singularities which is either diffeomorphic to a spherical 4-orbifold, or  covered by $\varepsilon_0$-caps and/ or   $\varepsilon_0$-necks; in the latter case  it must appear in the list in the conclusion of Proposition 6.2. To recover our original manifold $X$ from these components, denoted by  $X_1, X_2, \cdot\cdot\cdot$, we must invert the surgery procedure, that is, do orbifold connected sums among these components. As in the proof of Main Theorem in [CTZ],
we can divide the orbifold connected sum procedure into two steps. The first step
is to resolve (via genuine orbifold connected sums which are not manifold connected sums) all orbifold
singularities of $X_1, X_2, \cdot\cdot\cdot$ which are introduced
pairwise during the surgery process, and get
smooth manifolds, denoted by $Y_1, Y_2, \cdot\cdot\cdot$. Using Lemma 5.2 in [CTZ]
and a theorem in [Mc], which says that any diffeomorphism of a
3-dimensional spherical space form is isotopic to an isometry,  as in [CTZ] one can show that each  $Y_i$ is
diffeomorphic to $\mathbb{S}^4$, $\mathbb{RP}^4$, or $\mathbb{S}^3
\times \mathbb{R} /G$, where $G$ is a  discrete
subgroup of the isometry group of the round cylinder
$\mathbb{S}^3\times \mathbb{R}$ on which $G$ acts freely. The second step is to perform
manifold connected sums on $Y_1, Y_2, \cdot\cdot\cdot$
to invert the part of  surgery which does not introduce orbifold
singularities. The conclusion is that $X$ is diffeomorphic to a possibly infinite connected sum of  $\mathbb{S}^4$, $\mathbb{RP}^4$ and manifolds of the form $\mathbb{S}^3 \times \mathbb{R} /G$.  Finally the finiteness of diffeomorphism types of the factors $\mathbb{S}^3 \times \mathbb{R} /G$ follows from a simple analysis similar to that in case i in  the proof of Theorem 1.3  and [Mc]. Thus the proof of Theorem 1.1 is completed.
\hfill{$\Box$}

\hspace *{0.4cm}

Finally we state a theorem which generalizes Corollary 5.3 in [CTZ] to the noncompact case.

\begin{thm} \label{thm 6.3}\ \
 Let $(\mathcal{O}, g_0)$ be a complete, connected 4-orbifold with at most isolated singularities, with uniformly
positive isotropic curvature and with bounded geometry. Then there
is a finite collection $\mathcal{F}$ of spherical 4-orbifolds  $\mathbb{S}^4//\Gamma$ with at most isolated singularities and orbifolds with at most isolated singularities of the form
$\mathbb{S}^3 \times \mathbb{R} //G$, where $G$ is a  discrete subgroup of the isometry group of the standard metric
on $\mathbb{S}^3\times \mathbb{R}$,  such that $\mathcal{O}$ is diffeomorphic
to a (possibly infinite) orbifold connected sum of members of $\mathcal{F}$, where all the orbifold connected sum occurs at smooth points.
\end{thm}

\noindent {\bf  Proof} \ \  The proof is a slight modification of that of Theorem 1.1; cf. the proof of  Corollary 5.3 in [CTZ].
\hfill{$\Box$}

\hspace *{0.4cm}

\noindent {\bf Acknowledgements}  I would like to thank Prof. Xi-Ping Zhu for
his encouragement. His suggestion on Claim 2 in the proof of Proposition 3.6 in [Hu1] is very important for me, as already acknowledged in [Hu1]. I would also like to thank the referees for their helpful comments. I'm partially supported by NSFC no.11171025  and by Laboratory of Mathematics and Complex Systems, Ministry of Education.

\hspace *{0.4cm}

\bibliographystyle{amsplain}

\noindent {\bf References}

\hspace *{0.1cm}

[A] T. Aubin, Some nonlinear problems in Riemannian
geometry, Springer, 1998.

[BBB$^+$] L. Bessi$\grave{e}$res, G. Besson, M. Boileau, S.
Maillot and J. Porti, Geometrisation of 3-manifolds, Europ. Math.
Soc. 2010.

[BBM] L. Bessi$\grave{e}$res, G. Besson and S. Maillot, Ricci flow on
 open 3-manifolds and positive scalar curvature,  Geometry and Topology  15 (2011), 927-975.

[BMP] M. Boileau, S. Maillot, and J. Porti, Three-dimensional orbifolds and their geometric
structures, Soci$\acute{e}$t$\acute{e}$ Math$\acute{e}$matique de
France, Paris, 2003.

[BW] C. B$\ddot{o}$hm, B. Wilking, Manifolds with positive curvature operators are space forms, Ann. Math. 167 (2008), 1079-1097.

[B] J. Borzellino, Riemannian geometry of orbifolds,
Ph. D. thesis,  UCLA, 1992.

[B93] J. Borzellino, Orbifolds of maximal diameter, Indiana Univ. Math. J. 42 (1993),  37-53.

[BB]  J. Borzellino,  V. Brunsden,  Elementary orbifold differential
topology, Topology Appl. 159 (2012), no. 17, 3583-3589.

[BZ] J. Borzellino and S.-H. Zhu, The splitting theorem for orbifolds, Illinois J. Math. 38 (1994), no.4, 679-691.

[BJ]  Th. Br$\ddot{o}$cker and K. J$\ddot{a}$nich, Introduction to Differential Topology, Cambridge University Press, 1982.

[CaZ] H.-D. Cao, X.-P. Zhu, A complete proof of the
Poincar$\acute{e}$ and geometrization conjectures- application of
the Hamilton-Perelman theory of the Ricci flow, Asian J. Math. 10
(2006), 165-492.

[CE]  J. Cheeger, D. Ebin, Comparison Theorems in
Riemannian Geometry, North-Holland (1975).

[CG] J. Cheeger, D. Gromoll, On the structure of complete manifolds of nonnegative curvature, Ann. Math. 96 (1972), 413-443.

[CTZ] B.-L. Chen, S.-H. Tang and X.-P. Zhu,  Complete classification of
compact four-manifolds with positive isotropic curvature,  J. Diff.
Geom. 91 (2012), 41-80.

[CZ1] B.-L. Chen, X.-P. Zhu, Uniqueness of the Ricci flow on complete noncompact manifolds, J. Diff. Geom.
74 (2006), 119-154.

[CZ2] B.-L. Chen, X.-P. Zhu, Ricci flow with surgery
on  four-manifolds with positive isotropic curvature, J. Diff. Geom.
74 (2006), 177-264.

[C]  G. Chen, Calculus on orbifolds, Journal of Sichuan University (Natural Science Edition) 41 (2004), no.5, 931-939. (In Chinese)

[Ch] Y.-J. Chiang, Harmonic maps of $V$ -manifolds, Ann. Global Anal. Geom. 8 (1990), no.3, 315-344.

[CK] B. Chow, D. Knopf, The Ricci flow: an introduction, Amer. Math. Soc., 2004.

[CuZ] S. Cuccagna, B. Zimmermann, On the mapping class group of spherical 3-orbifolds, Proc. Amer. Math. Soc. 116 (1992), no.2, 561-566.

[D] D. DeTurck, Deforming metrics in the direction of their Ricci tensors, J. Diff. Geom. 18 (1983), no.1, 157-162; an improved version appeared later in Collected papers on Ricci flow, 163-165, International Press 2003, edited by H.D. Cao, B. Chow, S. C. Chu and S.T. Yau.

[F] C. Farsi, Orbifold spectral theory, Rocky
Mountain J. Math. 31 (2001), no.1, 215-235.

[GM] D. Gromoll, W. Meyer, On compact open manifolds of positive curvature, Ann. Math. 90 (1969), 75-90.

[GW] D. Gromoll, G. Walschap, Metric foliations and curvature, Progress in Mathematics 268, Birkh$\ddot{a}$user 2009.

[H82] R. Hamilton, Three-manifolds with positive
Ricci curvature, J. Diff. Geom. 17(1982), 255-306.

[H86] R. Hamilton, Four-manifolds with positive
curvature operator, J. Diff. Geom. 24(1986), 153-179.

[H93] R. Hamilton, The Harnack estimate for the Ricci flow, J. Diff. Geom. 37 (1993), 225-243.

[H95] R. Hamilton, A compactness property for solutions of the Ricci flow, Amer. J.
Math. 117 (1995), no.3, 545-572.

[H97] R. Hamilton, Four-manifolds with positive isotropic curvature, Comm. Anal. Geom. 5 (1997), 1-92; also in Collected Papers on Ricci flow, 342-407, edited by H. D. Cao, B. Chow, S. C. Chu and S. T. Yau, International Press 2003.

[H03] R. Hamilton, Three-orbifolds with positive Ricci curvature, in  Collected papers on Ricci flow,  521-524, edited by H.D. Cao, B. Chow, S. C. Chu and S.T. Yau, International Press 2003.

[He] E. Hebey,
Nonlinear analysis on manifolds: Sobolev spaces and inequalities,
Courant Institute of Mathematical Sciences, 1999.

[Hu1] H. Huang, Ricci flow on open 4-manifolds with positive isotropic
curvature, J. Geom. Anal. 23 (2013), no.3, 1213-1235.

 [Hu2] H. Huang, Three-orbifolds with positive scalar curvature, arXiv:1210.7331.

[HS]  G. Huisken, C. Sinestrari, Mean curvature flow with surgeries of two-convex hypersurfaces, Invent. Math. 175 (2009), no.1, 137-221.

[KL1] B. Kleiner, J. Lott, Notes on Perelman's
papers, Geom. Topol. 12 (2008), 2587-2855;  arXiv:0605667v5.

[KL2]  B. Kleiner, J. Lott, Geometrization of three-dimensional orbifolds via Ricci flow, arXiv:1101.3733v3, to appear in Ast$\acute{e}$risque.

[K] A. Kosinski, Differential Manifolds, Academic Press, 1993.

[Ko]  B. Kotschwar, An energy approach to the problem of uniqueness for the Ricci
flow, Comm. Anal. Geom. 22 (2014), no.1, 149-176.

[LSU] O. A. Lady$\check{z}$enskaja, V. A. Solonnikov, N. N. Ural'ceva, Linear and quasilinear equations of parabolic type, Amer. Math. Soc., 1968.

[L]  J. M. Lee, Introduction to smooth manifolds,
GTM 218, Springer, 2003.

[Lu] P. Lu, A compactness property for solutions of
the Ricci flow on orbifolds, Amer. J. Math.  123 (2001), 1103-1134.

[LT] P. Lu, G. Tian, Uniqueness of standard
solutions in the work of Perelman,
http://math.berkeley.edu/~lott/ricciflow/StanUniqWork2.pdf

[M] W. S. Massey, Algebraic topology: an introduction, GTM 56, Springer 1977.

[Mc] D. McCullough, Isometries of elliptic 3-manifolds, J. London Math. Soc. (2) 65 (2002), no. 1, 167--182.

[MM] M. Micallef, J.D. Moore, Minimal two-spheres
and the topology of manifolds with positive curvature on totally
isotropic two-planes, Ann. Math. (2) 127 (1988), 199-227.

[MW] M. Micallef, M. Wang, Metrics with nonnegative isotropic curvature, Duke. Math. J. 72 (1993), no. 3.  649-672.

[MT] J. Morgan, G. Tian, Ricci flow and the
Poincar$\acute{e}$ conjecture, Clay Mathematics Monographs 3, Amer.
Math. Soc., 2007.

[N] A. Naber, Noncompact shrinking four solitons with nonnegative curvature, J. Reine Angew. Math. 645 (2010), 125-153.

[Na] Y. Nakagawa, An isoperimetric inequality for
orbifolds, Osaka J. Math. 30 (1993), 733-739.

[P1] G. Perelman, The entropy formula for the Ricci flow and its geometric applications,
arXiv:math.DG/0211159.

[P2] G. Perelman, Ricci flow with surgery on three-manifolds, arXiv:math.DG/0303109.

[R] J. Ratcliffe, Foundations of hyperbolic manifolds, Second edition, GTM 149, Springer, 2006.

[S] W.-X. Shi, Deforming the metric on complete
Riemannian manifolds, J. Diff. Geom. 30 (1989), 223-301.

[T] W. Thurston, The geometry and topology of 3-manifolds, Princeton lecture notes (1979).

[To] P. Topping, Remarks on Hamilton's compactness theorem for Ricci flow, J. Reine Angew. Math. 692 (2014), 173-191.

 [Z]  Q. S. Zhang, Sobolev inequalities, heat
kernels under Ricci flow, and the Poincar$\acute{e}$ conjecture, CRC
Press 2011.

\vspace *{0.4cm}

School of Mathematical Sciences, Beijing Normal University,

Laboratory of Mathematics and Complex Systems, Ministry of Education,

Beijing 100875, P.R. China

 E-mail address: hhuang@bnu.edu.cn

\end{document}